\documentclass[11pt,a4paper]{article}
\usepackage[utf8]{inputenc}

\usepackage{graphicx}
\graphicspath{ {./images/} }

\title{Concentration of the maximum size of an induced subtree in moderately sparse random graphs}
\author{
  J.C. Buitrago Oropeza\thanks{Moscow Institute of Physics and Technology (National Research University), Department of Applied Mathematics and Informatics. E-mail address: buitrago.okh@phystech.edu}\\
}
\date{}

\usepackage{amsmath,amsthm,amssymb}
\usepackage{amsfonts}
\usepackage{mathrsfs}
\usepackage{mathtext}
\usepackage{mathtools}

\usepackage[left=2.5cm, right=2.5cm, top=2cm]{geometry}
\usepackage[T1]{fontenc}
\usepackage[english]{babel}
\usepackage{xcolor}
\usepackage[colorlinks,allcolors=blue]{hyperref}

\newtheorem{theorem}{Theorem}[section]

\theoremstyle{remark}

\newcommand{\RNumb}[1]{\uppercase\expandafter{\romannumeral #1\relax}}

\usepackage{color}   
\usepackage{hyperref}
\hypersetup{
    colorlinks=true, 
    linktoc=all,     
    linkcolor=blue,  
}
\usepackage{comment}

\usepackage{titling}
\thanksmarkseries{arabic}

\begin{document}
\maketitle

\begin{abstract}
Kamaldinov, Skorkin, and Zhukovskii proved that the maximum size of an induced subtree in the binomial random graph $G(n,p)$ is concentrated at two consecutive points, whenever $p\in(0,1)$ is a constant. Using improved bounds on the second moment of the number of induced subtrees, we show that the same result holds when $n^{-\frac{e-2}{3e-2}+\varepsilon}\leq p=o(1)$.
\end{abstract}

{\bf Keywords:} binomial random graph, maximum subgraph, concentration.

\section{Introduction}
In this paper, we discuss several properties of the binomial random graph $G(n,p)$ (see~\cite{1},\cite{2}, \cite{3}), where $p = p(n)$ is a function of $n$. Consider the set of vertices $[n] := \{1, \ldots, n\}$. In the random graph $G(n,p)$ on this set of vertices, each of the $\binom{n}{2}$ pairs of vertices is joined by an edge independently with probability $p$. More formally, $G(n,p)$ is a random element with values in the set of all graphs on $[n]$, following the distribution
\[
\mathbb{P}(G(n,p) = H) = p^{e(H)}(1-p)^{{n\choose 2} - e(H)},
\]
where $e(H)$ denotes the number of edges in $H$. 

A subset $U \subset [n]$ is called an \emph{independent set} if it contains no pair of adjacent vertices, i.e., the induced subgraph on $U$ contains no edges. The size of the largest independent set is called the \emph{independence number} of the graph. The distribution of the independence number in random graphs has been extensively studied. In~\cite{7,8}, it was shown that \textit{asymptotically almost surely} (a.a.s.), that is, with probability tending to 1 as $n\to \infty$, the independence number of $G(n,p = \mathrm{const})$ takes one of two values, $f_0(n)$ or $f_0(n)+1$, where
\[
f_0(n) = \left\lfloor 2\log_b n - 2\log_b\log_b n + 2\log_b\frac{e}{2} + 0.9 \right\rfloor, \quad b = \frac{1}{1-p}.
\]
Refinements of this result can be found in~\cite{9}. In such cases, we say that the independence number is \emph{concentrated at two points}. 

Later developments extended the two-point concentration phenomenon to other properties of random graphs. A natural generalization of the independence number problem is as follows: what is the maximum size of an \emph{induced} subgraph in $G(n,p)$ satisfying a given property?\footnote{A subset $A \subset V(G)$ of vertices induces in $G$ a subgraph $G|_A$, whose vertex set is $A$ and whose edge set consists of all edges of $G$ with both endpoints in $A$.} In particular, does the two-point concentration hold under a weaker constraint on the number of edges in the induced subgraph?\\

Let $t \colon \mathbb{Z}_{\ge 0} \to \mathbb{Z}_{\ge 0}$ be a function. Define $X_n[t]$ to be the maximum integer $k$ such that $G(n,p)$ contains an induced subgraph on $k$ vertices with {\it at most} $t(k)$ edges, and let $Y_n[t]$ be the maximum $k$ such that $G(n,p)$ contains an induced subgraph on $k$ vertices with {\it exactly} $t(k)$ edges. In~\cite{10}, it was shown that under certain conditions on $t$, the value $X_n[t]$ is concentrated at two points.
\begin{theorem}[N.~Fountoulakis, R.~J.~Kang, C.~McDiarmid, 2014]
Let $t = t(k) = o\left( \frac{k \sqrt{\ln k}}{\sqrt{\ln\ln k}} \right)$ and $p = \mathrm{const} \in (0,1)$. Define
\[
f_t(n) = \left\lfloor 2\log_b n + (t - 2)\log_b\log_b (np) - t\log_b t + t\log_b\frac{2pe}{1-p} + 2\log_b\frac{e}{2} + 0.9 \right\rfloor.
\]
Then $X_n[t] \in \{f_t(n), f_t(n)+1\}$ a.a.s.
\label{McD}
\end{theorem}

On the other hand, in~\cite{11}, it was shown that the random variable $Y_n[t]$ does not exhibit concentration on a finite set of points when $t(k)$ is close to $p\binom{k}{2}$.

\begin{theorem}[J.~Balogh, M.~Zhukovskii, 2019]
Let $t(k) = p\binom{k}{2} + O(k)$.
\begin{enumerate}
\item[(i)] There exists a constant $\mu > 0$ such that for any $c > \mu$ and $C > 2c + \mu$, we have
\begin{multline*}
0 < \liminf_{n \to \infty} \mathbb{P}\left( n - C \sqrt{\frac{n}{\ln n}} < Y_n[t] < n - c \sqrt{\frac{n}{\ln n}} \right)\\
\le \limsup_{n \to \infty} \mathbb{P}\left( n - C \sqrt{\frac{n}{\ln n}} < Y_n[t] < n - c \sqrt{\frac{n}{\ln n}} \right) < 1.
\end{multline*}
\item[(ii)] Suppose that for any sequence $m_k = O\left( \sqrt{k / \ln k} \right)$ of nonnegative integers, the following holds:
\[
\left| \left( t(k) - p\binom{k}{2} \right) - \left( t(k - m_k) - p\binom{k - m_k}{2} \right) \right| = o(k).
\]
Then for any $\varepsilon > 0$, there exist constants $c,C > 0$ such that
\[
\liminf_{n \to \infty} \mathbb{P}\left( n - C \sqrt{\frac{n}{\ln n}} < Y_n[t] < n - c \sqrt{\frac{n}{\ln n}} \right) > 1 - \varepsilon.\\
\]
\end{enumerate}
\label{large_edges}
\end{theorem}

Finally, in~\cite{12}, the two-point concentration for $Y_n[t]$ was established for small values of $t$.

\begin{theorem}[D.~Kamaldinov, A.~Skorkin, M.~Zhukovskii, 2021]
Let $R > 0$. Then there exists $\varepsilon > 0$ such that for any sequence of nonnegative integers $t(k)$ satisfying the following conditions:
\begin{itemize}
\item $|t(k+1)/t(k) - 1| \le \frac{R}{k}$ for all sufficiently large $k$,
\item $t(k) < \varepsilon k^2$ for all sufficiently large $k$,
\end{itemize}
there exists a function $f(n)$ satisfying $|f(n) - 2\log_{1/(1-p)} n| \le (3\varepsilon \ln(1/\varepsilon))\ln n$, such that a.a.s. $Y_n[t] \in \{f(n), f(n)+1\}$.
\end{theorem}

The question of the two-point concentration phenomenon has also been raised in the context of structural constraints on the induced subgraph. For example, in~\cite{12}, the maximum size of an induced tree in $G(n, p = \mathrm{const})$ was investigated and its two-point concentration was established.

\begin{theorem}[D.~Kamaldinov, A.~Skorkin, M.~Zhukovskii, 2021]
\label{thm1}
Let 
$$
f_{\varepsilon}(n) = \lfloor 2\log_{1/(1-p)} (enp) + 2 + \varepsilon \rfloor.
$$
Then there exists $\varepsilon > 0$ such that a.a.s. the maximum size of an induced tree in $G(n, p = \mathrm{const})$ takes one of the two values: $f_{\varepsilon}(n)$ or $f_{\varepsilon}(n)+1$.
\end{theorem}

In~\cite{13}, a similar result was obtained for the maximum size of an induced forest.

\begin{theorem}[M.~Krivoshapko, M.~Zhukovskii, 2021]
There exists $\varepsilon > 0$ such that a.a.s. the maximum size of an induced forest in $G(n, p = \mathrm{const})$ takes one of the two values:
\[
\lfloor 2\log_{1/(1-p)} (enp) + 2 + \varepsilon \rfloor
\quad \text{or} \quad
\lfloor 2\log_{1/(1-p)} (enp) + 3 + \varepsilon \rfloor.
\]
\end{theorem}

Let us notice that the above results concern random graphs with constant edge probability. This naturally leads to the question of whether similar results hold in the sparse regime, where $p(n) = o(1)$. In~\cite{14}, it was shown that the independence number of $G(n,p)$ is concentrated at two points when $n^{-2/3 + \varepsilon} < p < 1/\log^2 n$.
\begin{theorem}[T.~Bohman, J.~Hofstad, 2023]\label{thm_bh}
Let $\varepsilon > 0$ and suppose $n^{-2/3 + \varepsilon} < p < 1/\log^2 n$. Then a.a.s. the independence number of the random graph $G(n,p)$ is concentrated at two points.
\end{theorem}

Moreover, it was shown that the bound on $p$ in Theorem \ref{thm_bh} is optimal, as the two-point concentration of the independence number fails to hold when $p = n^{\gamma}$, for $\gamma \leq -2/3$.

The study of the independence number has also been extended to other models. In \cite{4}, a lower bound was established on the independence number of the random subgraph $G_{1/2}(n,2,1)$, where $G_{p}(n,2,1)$ is obtained by selecting each edge of the distance graph $G(n,2,1)$ independently with probability~$p$.

\paragraph{Related work on maximum induced trees.} The study of maximum induced trees in $G(n,p)$ dates back to 1983, when Erd\H{o}s and Palka \cite{EZ} showed that, for constant $p$, the maximum size of an induced tree is asymptotically $2(1+o(1))\log_{1/(1-p)}(np)$. Shortly after, Fernández de la Vega \cite{Vega} and Frieze and Jackson \cite{FJ} extended this result to the sparse regime $p = c/n$ with $c > 1$, showing that the size of such a subtree is linear in $n$. Later, in \cite{Vega2}, Fernández de la Vega refined this result --- in particular, from his argument the asymptotics of the maximum size of an induced tree follows for all $1/n\ll p\leq 1-\varepsilon$. More recently, Dutta and Subramanian \cite{DS} revisited this line of investigation and improved the error term for all $n^{-1/2}(\ln n)^2 \leq p \leq 1 - \varepsilon$. They showed that a.a.s. the maximum size of an induced subtree in $G(n, p)$ is $2\log_{1/(1-p)}(np) + O(1/p)$. 

In 2020, Dragani\'{c} \cite{Drag} proved that, for any given tree $T$ of bounded maximum degree and size at most $(1-\varepsilon)\log_{1/(1-p)}(np)$, and for $n^{-1/2}(\ln n)^2 \leq p \leq 0.99$, the random graph $G(n,p)$ a.a.s. contains an induced copy of $T$. Finally, Akhmejanova, Kozhevnikov, and Zhukovskii~\cite{AKZ} in 2024 proved 2-point concentration of the maximum size of an induced tree with maximum degree at most $\Delta$ in dense binomial random graphs $ G(n, p)$ with constant probability $p$.

We also refer an interested reader to similar results on maximum sizes of induced matchings, cycles, and forests~\cite{AK,AKZ,16,CDKS,DGK,DS,Glock}.



\paragraph{Main contributions of the paper.} Extending this line of research, we prove in this paper that Theorem~\ref{thm1}, which concerns the maximum size of an induced tree in $G(n, p)$, holds for a fairly broad range of values of $p = o(1)$.

\begin{theorem}
\label{th:new}
Let $\varepsilon>0$, $n^{-\frac{e-2}{3e-2}+\varepsilon}\leq p=o(1)$. Then there exists a constant $\delta>0$ such that a.a.s. the maximum size of an induced tree in $G(n,p)$ belongs to $\{g(n),g(n)+1\}$, where
\begin{equation*}
    g(n)= \lfloor 2\log_{1/(1-p)}(enp)+ \delta \rfloor.
\end{equation*}
\end{theorem}

We largely follow the proof framework developed in~\cite{12}, where the two-point concentration for the size of maximum induced subtrees was established in the case $p = \text{const}$ (see Theorem~\ref{thm1}). However, in the sparse regime, where $p = o(1)$, the asymptotic behavior of the second moment of the number of induced tree copies changes significantly, necessitating improved bounds. The main difficulty in deriving a sufficient bound on the second moment is in estimating the number of trees on $k=\Theta((\ln n)/p)$ vertices containing a given forest of size $\ell$. In particular, bounds from~\cite{12} are not enough in the case when the given forest contains almost all vertices of the tree, namely $\ell\geq k-O((\ln k)/p)$.

It is worth noting that we suspect that the bound on $p$ established here is not optimal, and that the two-point concentration may, in fact, hold for smaller values of $p$ as well. However, proving this would likely require new ideas, as the obtained requirement on $p$ is necessary for our bounds on the number of overlapping trees to imply that the variance of the number of induced trees is sufficiently small.

\paragraph{Structure of the paper.} The remainder of the paper is dedicated to the proof of Theorem~\ref{th:new}. The proof is structured as follows. 

In Section~\ref{comp_moment}, we define $X_k$ as the number of induced subtrees of size~$k$ in the random graph $G(n,p)$. We identify a value $\hat{k}$ such that 
$$
{\sf E}X_{\lceil \hat{k} + \delta \rceil}\to 0
\quad\text{ and }\quad
{\sf E}X_{\lfloor \hat{k} - 1 + \delta \rfloor}\to\infty\quad\text{ as }n\to\infty.
$$
Due to Markov's inequality, the first bound immediately implies that ${\sf P}(X_{\lceil \hat{k} + \delta \rceil} \geq 1)\to 0$ as $n\to\infty$, as needed. In order to derive that ${\sf P}(X_{\lfloor \hat{k} - 1 + \delta \rfloor}\geq 1)\to 1$, we exploit the usual bound
\begin{equation}\label{the_first}
 \frac{\mathrm{Var}X_k}{({\sf E}X_k)^2} \leq \sum_{\ell = 2}^{k-1} \frac{F_{\ell}}{({\sf E}X_k)^2},
\end{equation}
where $F_{\ell}$ is the number of pairs of $k$-trees on $[n]$ that share exactly $\ell$ vertices.

 In subsequent sections, we derive suitable upper bounds for $F_{\ell}/({\sf E}X_k)^2$ in order to show that the right-hand side of~\eqref{the_first} is $o(1)$.

In particular, in Section~\ref{p_small} we address the case $n^{-\frac{e-2}{3e-2}+\varepsilon} \leq p < \frac{1}{2\ln n}$. In this regime, we partition the set $\{2, \ldots, k-1\}$ of possible values of $\ell$ into four subsets and verify that the corresponding contribution to the sum in~\eqref{the_first} is $o(1)$, in each subset separately. Finally, Section~\ref{p_large} considers the case $\frac{1}{2\ln n} \leq p = o(1)$ and likewise establishes that the right-hand side of~\eqref{the_first} is~$o(1)$.\\

A short preliminary version of this paper appeared without proofs in~\cite{17}.

\section{Computation of moments}
\label{comp_moment}

Let $X_{k}$ be the number of induced subtrees in $G(n,p)$ of size $k$.  We have that
\begin{equation}
    {\sf E}X_k = \binom{n}{k} k^{k-2} p^{k-1}(1-p)^{\binom{k}{2}-k+1}.
    \label{eq:expectation_X_k}
\end{equation}

Note that $g(n)$ from the statement of Theorem~\ref{th:new} is asymptotically smaller than $\sqrt{n}$. Since a graph without induced trees of size $k$ does not contain induced trees of larger size, throughout the proof, we can restrict ourselves with subtrees of size $o(\sqrt{n})$. Also, the right-hand side in~\eqref{eq:expectation_X_k} approaches infinity when $k=O(1)$. We will show that there exists a minimum $k=o(\sqrt{n})$ such that $\mathbb{E}X_k\leq 1$, and $\mathbb{E}X_k=o(1)$ for all larger~$k=o(\sqrt{n})$. 

Consider $k$ such that $1\ll k\ll \sqrt{n}$. We will use the following well-known asymptotic relation that follows immediately from Stirling's formula:  $\binom{n}{k}\sim \frac{1}{\sqrt{2\pi k}} \left( \frac{ne}{k}\right)^k$. For the sake of completeness, let us give a short proof. Indeed, since 
$$
\left( \frac{n}{n-k}\right)^{n-k}=e^{(n-k)\ln (1+\frac{k}{n-k})}=
e^{ k+O\left(\frac{k^2}{n-k}\right)}=e^{k+o(1)},
$$
then, applying Stirling's approximation, we get
\begin{align*}
    \binom{n}{k} &\sim \frac{\sqrt{2\pi n} (n/e)^n}{\sqrt{2\pi k} (k/e)^k \sqrt{2\pi (n-k)} ((n-k)/e)^{(n-k)}}\\
    &= \frac{1}{\sqrt{2 \pi k}}\sqrt{\frac{n}{n-k}} \left(\frac{n}{k} \right)^{k} \left( \frac{n}{n-k}\right)^{n-k} \sim \frac{1}{\sqrt{2\pi k}} \left( \frac{ne}{k}\right)^k.
\end{align*}
We thus get ${\sf E}X_k \sim e^{\gamma(k)}$, where
$$
\gamma(k)= -\frac{1}{2}\ln 2\pi +k\ln n + k - \frac{5}{2} \ln k + (k-1) \ln \left(\frac{p}{1-p}\right) - \binom{k}{2} \ln \left(\frac{1}{1-p} \right).
$$
Notice that $\frac{\partial \gamma(k)}{\partial k}<0$, if $k>\frac{\ln n + \ln p+1+o(1)}{\ln [1/(1-p)]}$, since
\begin{equation*}
    \begin{split}
        \frac{\partial \gamma(k)}{\partial k} & = \ln n +1 - \frac{5}{2k} + \ln \left( \frac{p}{1-p}\right) - \left(k-\frac{1}{2} \right) \ln  \left( \frac{1}{1-p}\right)\\
        & =\ln n +\ln p -k \ln \left(\frac{1}{1-p} \right)+1- \frac{5}{2k}+ \ln \left( \frac{1}{(1-p)^{3/2}}\right) \\
        & = \ln n + \ln p -k \ln\left(\frac{1}{1-p} \right)+1+o(1).
    \end{split}
\end{equation*}
So, $\frac{\partial \gamma(k)}{\partial k}$ changes its sign from positive to negative at $ \frac{\ln n + \ln p+1+o(1)}{\ln [1/(1-p)]}$. Since $\gamma(k)\to\infty$ for small $k$, it means that either $\gamma$ never reaches 0, or it reaches $0$ at some point, and after that becomes negative. In other words, if $\gamma(k)=0$ has a solution, then this solution is unique.
In the next section, we show that such a solution, that we will denote by $\hat{k}$, exists.

\subsection{Computation of $\hat{k}$}
Let us rewrite the expression for $\gamma(k)$ in a more convenient way:
\begin{equation*}
    \begin{split}
        \gamma(k) & = -\frac{1}{2}\ln 2\pi +k\ln n + k - \frac{5}{2} \ln k + (k-1) \ln \left(\frac{p}{1-p}\right) - \binom{k}{2} \ln \left(\frac{1}{1-p} \right) \\
        & = k \left[\ln n +1 + \ln \left( \frac{p}{1-p} \right)+\frac{1}{2}\left( 1-k\right) \ln \left(\frac{1}{1-p} \right)\right] -\frac{1}{2}\ln 2\pi- \frac{5}{2} \ln k - \ln \left( \frac{p}{1-p}\right)\\
        & =k \left[\ln(np) +1 + \frac{3}{2}\ln \left( \frac{1}{1-p} \right)-\frac{k}{2} \ln \left(\frac{1}{1-p} \right)\right] -\frac{1}{2}\ln 2\pi- \frac{5}{2} \ln k - \ln \left( \frac{p}{1-p}\right).
    \end{split}
\end{equation*}
Letting 
$$
k^*=\frac{2}{\ln[1/(1-p)]} \left[ \ln (np) +1 +\frac{3}{2}\ln \left( \frac{1}{1-p} \right)\right],
$$
we get 
$$
\gamma(k^*)=-\frac{1}{2}\ln 2\pi- \frac{5}{2} \ln k^* - \ln \left( \frac{p}{1-p}\right)<0.
$$
Thus, following our conclusions on monotonicity of $\gamma$, we get that the desired unique solution $\hat{k}=o(\sqrt{n})$ of the equation $\gamma=0$ exists. We may now find $\varepsilon=\varepsilon(n) >0$ such that $\varepsilon=O(1)=o(k^*)$ and $\gamma(k^*-\varepsilon)=0$. For any $\varepsilon=O(1)$ (it can be negative),
\begin{equation}
    \begin{split}
    \gamma(k^*- \varepsilon) & =  (k^*-\varepsilon)\left(\frac{\varepsilon}{2} \ln \left( \frac{1}{1-p} \right) \right) - \frac{1}{2} \ln 2\pi - \frac{5}{2} \ln (k^*- \varepsilon) - \ln \left(\frac{p}{1-p}\right)\\
    &= \frac{k^* \varepsilon p}{2}-\frac{5}{2}\ln k^*-\ln p -\frac{1}{2}\ln 2\pi+ o(1).
    \end{split}
\label{eq:gamma_change}
\end{equation}
Thus, $\gamma(\hat k)=0$ where $\hat k=k^*-\varepsilon$ and 
\begin{equation*}
    \begin{split}
        \varepsilon & =\frac{2}{k^*p}\left[\frac{5}{2}\ln k^*+\ln p +\frac{1}{2}\ln 2\pi + o(1) \right]\\
        & =\frac{1}{k^*p} \left[5\ln \left( 2\frac{\ln np}{p}\right)+2\ln p +\ln 2\pi +o(1)\right]\\
        & =\frac{1}{k^*p} \left[-3\ln p +5\ln \ln np +5\ln 2 + \ln 2\pi + o(1)\right]\\
        & = -\frac{3\ln p}{k^*p}+o(1)
         = -\frac{3\ln p}{2 \ln( np)}+o(1).
    \end{split}
\end{equation*}
 We finally conclude that
\begin{equation}
    \begin{split}
       \hat{k} & = \frac{2}{\ln[1/(1-p)]} \left[ \ln (np) +1 +\frac{3}{2}\ln \left( \frac{1}{1-p} \right)\right]+\frac{3\ln p}{2 \ln(np)}+o(1)\\
        & =2\log_{1/(1-p)}(enp)+ \frac{3\ln p}{2\ln(np)} + 3 + o(1). 
    \end{split}  
\label{value_k}
\end{equation}

\subsection{First and second moment method}
\label{method_moment}
From~\eqref{eq:gamma_change}, it follows that, for any constant $\delta>0$,
\begin{equation}
{\sf E}X_{\lceil \hat{k}+\delta \rceil}\sim e^{\gamma(\lceil \hat{k}+\delta \rceil)} \leq e^{\gamma(\hat{k}+\delta)}
=e^{\gamma(k^*-(\varepsilon-\delta))}=(1+o(1))e^{\gamma(\hat k)-k^*\delta p/2}=
(1+o(1))e^{-k^*\delta p/2}=o(1).
\label{eq:expectation_upper_0}
\end{equation}
In a similar way, for any constant $\delta\in(0,1)$,
$$
{\sf E}X_{\lfloor \hat{k}-1 +\delta \rfloor}\sim e^{\gamma(\lfloor \hat{k}-1+\delta \rfloor)} \geq e^{\gamma(\hat{k}-1+\delta)}
=(1+o(1))e^{\gamma(\hat k)+k^*(1-\delta)p/2}=
(1+o(1))e^{k^*(1-\delta) p/2}\to\infty
$$
as $n\to\infty$.

First of all, by Markov's inequality,~\eqref{eq:expectation_upper_0} implies
$$
{\sf P}(X_{\lceil \hat{k}+\delta \rceil} \geq 1) \leq {\sf E}X_{\lceil \hat{k}+\delta \rceil} \stackrel{n\to\infty}\to 0.
$$
It remains to prove that 
\begin{equation}
    {\sf P}(X_{\lfloor \hat{k}-1 +\delta \rfloor} \geq 1) \to 1.
    \label{eq:Cheb_to_1}
\end{equation}

For simplicity of notation, set $k= \lfloor \hat{k}-1 +\delta \rfloor$. Due to Chebyshev's inequality, we have

\begin{equation}\label{eq_2.1}
{\sf P}(X_k=0)\leq \frac{\mathrm{Var}X_k}{({\sf E}X_k)^2} = \frac{{\sf E} X_k (X_k -1) - ({\sf E} X_k)^2}{({\sf E} X_k)^2} +\frac{1}{{\sf E}X_k} = \frac{{\sf E} X_k (X_k -1) - ({\sf E} X_k)^2}{({\sf E} X_k)^2}+o(1).
\end{equation}

Our goal is to show that (\ref{eq_2.1}) is $o(1)$. For that, it remains to find a suitable upper bound for ${\sf E} X_k (X_k -1)-({\sf E} X_k)^2 $.\\ 

In the usual way, we represent the random variable $X_k$ as a sum of indicators $I_{B(U)}$ over $U\in{[n]\choose k}$, where $B(U)$ is the event that the set of vertices $U$ induces a tree in $G(n,p)$. Therefore
\begin{align}
{\sf E} X_k (X_k -1) &= \sum_{U \neq U'}\mathbb{P}(B(U)\cap B(U'))\notag \\
&=\sum_{\ell\in\{0,1\}} {n\choose k}{k\choose \ell}{n-k\choose k-\ell}(k^{k-2})^2p^{2(k-1)}(1-p)^{2{k\choose 2}-2(k-1)}+\sum_{\ell=2}^{k-1}G_{\ell},
\label{ecuacuion_2}
\end{align}
where
$$
G_{\ell}=\sum_{r=0}^{\ell-1} {n\choose k}{k\choose \ell}{n-k\choose k-\ell} N(k,\ell,r) p^{2(k-1)-r}(1-p)^{2{k\choose 2}-{\ell\choose 2}-2(k-1)+r},
$$
and $N(k,\ell,r)$ is the number of pairs of labeled trees on $[k]$ and $[2k-\ell]\setminus[k-\ell]$ such that the number of edges in their intersection on $\ell$ vertices is exactly $r$. We can further bound $G_{\ell}$ in the following way for every $\ell\geq 2$:
\begin{equation}\label{eq_1}
    G_{\ell} \leq \displaystyle \binom{n}{k} \binom{k}{\ell} \binom{n-k}{k-\ell} \ell \displaystyle \max_{r \in \{0,\ldots,\ell-1 \}}p^{2(k-1)-r} (1-p)^{2 \binom{k}{2} - \binom{\ell}{2} -2(k-1)+r} N(k,\ell,r).
\end{equation}

Now, we also want to express $({\sf E} X_k)^2$ as a sum depending on $\ell$.

Replacing $\binom{n}{k}$ by $\displaystyle \sum_{\ell =0}^k \binom{k}{\ell} \binom{n-k}{k-\ell}$ in  $({\sf E} X_k)^2$ we obtain
\begin{equation}\label{ecuacion_3}
    \begin{split}
        ({\sf E} X_k)^2 & = \binom{n}{k} \displaystyle \sum_{\ell =0}^k \binom{k}{\ell} \binom{n-k}{k-\ell} k^{2(k-2)} p^{2(k-1)}(1-p)^{2\binom{k}{2}-2k+2}.
    \end{split}   
\end{equation}

Combining (\ref{ecuacuion_2}), (\ref{eq_1})  and  (\ref{ecuacion_3}), we thus get
\begin{equation*}\label{ecuacion_4}
\begin{split}
    & {\sf E} X_k (X_k -1) - ({\sf E} X_k)^2 \leq \\
    & \leq   \displaystyle \sum_{\ell =2}^{k-1} \binom{n}{k} \binom{k}{\ell} \binom{n-k}{k-\ell} p^{2(k-1)} (1-p)^{2 \binom{k}{2} -2(k-1)} \left[ \displaystyle\ell \max_{r}  \frac{(1-p)^{r-\binom{\ell}{2}}}{p^r} N(k, \ell, r) -k^{2(k-2)} \right]\\
     & \leq   \displaystyle \sum_{\ell =2}^{k-1} \binom{n}{k} \binom{k}{\ell} \binom{n-k}{k-\ell} p^{2(k-1)} (1-p)^{2 \binom{k}{2} - \binom{\ell}{2} -2(k-1)} \ell
     \max_{r}  \frac{(1-p)^{r}}{p^r}N(k, \ell, r).
\end{split}  
\end{equation*}
Letting 
$$
F_{\ell}:=\binom{n}{k} \binom{k}{\ell} \binom{n-k}{k-\ell} p^{2(k-1)} (1-p)^{2 \binom{k}{2} - \binom{\ell}{2} -2(k-1)} \ell \max_{r\in\{0,\ldots,\ell-1\}}  \frac{(1-p)^{r}}{p^r}N(k, \ell, r),
$$
we get
\begin{equation}\label{ecuacion_4_1}
\frac{{\sf E} X_k (X_k -1) - ({\sf E} X_k)^2}{({\sf E} X_k)^2} \leq \displaystyle \sum_{\ell =2}^{k-1} \frac{F_{\ell}}{({\sf E}X_k)^2}.
\end{equation}

Our objective is to establish that the right-hand side of (\ref{ecuacion_4_1}) is $o(1)$ for all $n^{-\frac{e-2}{3e-2}+\varepsilon} \leq p = o(1)$. We will first prove this for $n^{-\frac{e-2}{3e-2}+\varepsilon} \leq p < \frac{1}{2\ln n}$ and then for $\frac{1}{2\ln n} \leq p = o(1)$. For the range $\frac{1}{2\ln n} \leq p = o(1)$, we will show that the argument from \cite{12}, originally formulated for constant $p$, can be suitably adapted to this setting with only a few adjustments. This is done in Section~\ref{p_large}. Nevertheless, the case $n^{-\frac{e-2}{3e-2}+\varepsilon} \leq p < \frac{1}{2\ln n}$ requires a special treatment, since bounds on $F_{\ell}$ from~\cite{12} do not suffice in this case. We prove the desired bound on $F_{\ell}$ and thus derive~\eqref{eq:Cheb_to_1} in this case, in Section~\ref{p_small}.

\section{Proof of Theorem~\ref{th:new} for $n^{-\frac{e-2}{3e-2}+\varepsilon}\boldsymbol{ \leq p <} \frac{1}{2\ln n}$}\label{p_small}

In order to prove that the right-hand side of (\ref{ecuacion_4_1}) is small, we partition the set $\{2,\ldots,k-1\}$ of values of $\ell$ into four parts, and then individually show that in each of these parts the sum is $o(1)$. Let us describe this partition. We consider a sequence $w(n) \to \infty$ such that $w(n)= o\left(\sqrt{\ln n}\right)$. Also, set 
\begin{equation}\label{l^*}
    \ell^*= \frac{2\ln (np) - 2\ln (4ek)}{\ln[1/(1-p)] }.
\end{equation}
Then, we define the following partition: 
$$
\{2,\ldots,k-1\}=[2,\ell^*]\sqcup \left(\ell^*, \left \lfloor k-\frac{w}{p} \right\rfloor\right]\sqcup\left( \left \lfloor k-\frac{w}{p} \right\rfloor, k- \frac{1}{2p}\right]\sqcup\left(k-\frac{1}{2p}, k-1\right].
$$
Each of these parts is considered in its own subsequent subsection. We shall use different upper bounds for $N(k, \ell, r)$ in order to achieve the intended result. In particular, in Section~\ref{sc:part1} we consider the first interval $[2,\ell^*]$ and the obvious bound $N(k, \ell, r)\leq k^{2(k-2)}$. In Sections~\ref{sc:part2}~and~\ref{sc:part4}, we consider the second and last intervals respectively, and the bound $N(k, \ell, r)\leq k^{k-2}f(k,\ell,r)$,
where $f(k,\ell,r)$ is an upper bound (we will define its value in that sections), over all choices of a forest $F$ with $r$ edges on $[\ell]$, on the number of trees on $[k]$ that induce $F$ on $[\ell]$. Finally, in Section~\ref{sc:part3}, we consider the remaining interval $(k-w/p, k- 1/(2p)]$ and use the bound $N(k, \ell, r)\leq \varphi(\ell,r)[f(k,\ell,r)]^2$, where $\varphi(\ell,r)$ is the number of forests on $[\ell]$ with $r$ edges.

\subsection{Part 1: $\ell \in [2,\ell^*]$}
\label{sc:part1}
Here we use the bound $N(k,\ell,r)\leq (k^{k-2})^2$. Thus we obtain for every $\ell \in [2,\ell^*]$ the following bound for $F_{\ell}$:
\begin{equation*}
    \begin{split}
         F_{\ell} & \leq \binom{n}{k} \binom{k}{\ell} \binom{n-k}{k-\ell} p^{2(k-1)} (1-p)^{2 \binom{k}{2} - \binom{\ell}{2} -2(k-1)} (k^{k-2})^2 \ell \max_{r \in\{0,\ldots,\ell-1\}}  \frac{(1-p)^{r}}{p^r}\\
         & \leq  \binom{n}{k} \binom{k}{\ell} \binom{n-k}{k-\ell}  p^{2(k-1)} (1-p)^{2 \binom{k}{2} - \binom{\ell}{2} -2(k-1)} (k^{k-2})^2 k \left(\frac{1-p}{p}\right)^{\ell}.
    \end{split}
\end{equation*}
We thus get
\begin{equation*}
\begin{split}
    \displaystyle \sum_{\ell =2}^{\ell^*} \frac{F_{\ell}}{({\sf E} X_k)^2} & \leq \displaystyle \sum_{\ell =2}^{\ell^*} \frac{\binom{n}{k} \binom{k}{\ell} \binom{n-k}{k-\ell} p^{2(k-1)} (1-p)^{2 \binom{k}{2} - \binom{\ell}{2} -2(k-1)} (k^{k-2})^2 k \left(\frac{1-p}{p}\right)^{\ell}}{ \binom{n}{k}^2  k^{2(k-2)} p^{2(k-1)}(1-p)^{2\binom{k}{2}-2k+2}} \\ 
    & =  \displaystyle \sum_{\ell =2}^{\ell^*} \frac{\binom{k}{\ell} \binom{n-k}{k-\ell}}{\binom{n}{k}} k (1-p)^{-\binom{\ell}{2}}  \left(\frac{1-p}{p} \right)^{\ell}.
\end{split}  
\end{equation*}    
    
By Stirling's formula, we have that for some constant $c<1$, for all sufficiently large $n$, and all $\ell\in[2,\ell^*]$,
\begin{equation}
    \begin{split}
      \frac{\binom{k}{\ell} \binom{n-k}{k-\ell}}{\binom{n}{k}} & <\frac{\binom{k}{\ell} \binom{n}{k-\ell}}{\binom{n}{k}}=
\frac{(k!)^2(n-k)!}{\ell!((k-\ell)!)^2(n-k+\ell)!}\leq c\cdot\frac{k^{2k}(n-k)^{n-k}}{\ell^{\ell}(k-\ell)^{2(k-\ell)}(n-k+\ell)^{n-k+\ell}}\\
&=c\left(\frac{k^2}{\ell n} \right)^{\ell}
\left(\frac{k}{k-\ell}\right)^{2(k-\ell)}
\left(\frac{n}{n-k+\ell}\right)^{\ell}
\left(\frac{n-k}{n-k+\ell}\right)^{n-k}\\
&\leq c\left(\frac{k^2}{\ell n} \right)^{\ell}
\exp\left[2\ell+\frac{\ell(k-\ell)}{n-k+\ell}-\frac{\ell(n-k)}{n-k+\ell}\right]=
c\left(\frac{k^2}{\ell n} \right)^{\ell}
\exp\left[\ell+\frac{\ell k}{n-k+\ell}\right]\\
&\leq \left( \frac{k^2e}{\ell n} \right)^{\ell},  
    \end{split}
    \label{bound_11}
\end{equation}
where the last inequality is due to $\ell k\leq k^2=o(n)$. Therefore,
\begin{equation}
\begin{split}
     \displaystyle \sum_{\ell =2}^{\ell^*} \frac{F_{\ell}}{({\sf E} X_k)^2} & \leq \displaystyle \sum_{\ell =2}^{\ell^*}  k \left( \frac{k^2e (1-p)}{\ell n p} \right)^{\ell}(1-p)^{-\binom{\ell}{2}} \leq \displaystyle \sum_{\ell =2}^{\ell^*} k \left( \frac{k^2e (1-p)}{\ell n p} (1-p)^{-\frac{\ell}{2}}\right)^{\ell} = \displaystyle \sum_{\ell =2}^{\ell^*}  e^{g(\ell)},
\end{split}  
\label{bound_1}
\end{equation}
where 
$$
g(\ell)= \ln k+ \ell\left(1+2\ln k +(1-\ell/2)\ln (1-p)-\ln n- \ln \ell - \ln p\right).
$$

Let us now evaluate the behavior of $g(\ell)$. From~(\ref{value_k}) we get $\ln k=-\ln p +\ln \ln (np)+O(1)$, and therefore
\begin{equation*}
    \begin{split}
       \frac{\partial}{\partial \ell} g(\ell) =& 2 \ln k+(1-\ell) \ln (1-p)-\ln n-\ln \ell-\ln p\\
       = &3\ln (1/p)+2\ln \ln (np)-\ln \ell +\ell p-\ln n + O(1)+O(p\ln n).  
    \end{split}
\end{equation*}
It follows that 
\begin{equation*}
    \begin{split}
       \left.\frac{\partial}{\partial \ell} g(\ell)\right|_{\ell=2}=&3\ln (1/p)+2\ln \ln (np)-\ln n + O(1)+O(p\ln n)\\
       =&-\ln (np^3)+O(\ln \ln n)+O(p\ln n). 
    \end{split}
\end{equation*}
As soon as $p\geq n^{-1/3+\varepsilon}$ (note that our requirement on $p$ is even stronger), we get $\left.\frac{\partial}{\partial \ell} g(\ell)\right|_{\ell=2}<0$.

On the other hand, from~(\ref{l^*}) we get $\ln \ell^*=-\ln p+\ln \ln (np)+O(1)$, and therefore
\begin{equation*}
    \begin{split}
      \left.\frac{\partial}{\partial \ell} g(\ell)\right|_{\ell=\ell^*} & =3\ln (1/p)+2\ln \ln (np)-\ln \ell^*-\ln n + 2\ln (np)-2\ln (4ek) + O(1)+O(p\ln n)\\
      & =-2\ln(1/p)+\ln n + O(\ln \ln n) +O(p\ln n)\\
      & = \ln (np^2) + O(\ln \ln n)+O(p\ln n).
    \end{split}
\end{equation*}
As soon as $p\geq n^{-1/2+\varepsilon}$, we get that $\left.\frac{\partial}{\partial \ell} g(\ell)\right|_{\ell=\ell^*}>0$.

Furthermore, since $\frac{\partial^2}{\partial \ell^2} g(\ell)=-1/\ell-\ln(1-p)$, then $\frac{\partial^2}{\partial \ell^2} g(\ell)$ changes its sign from negative to positive at $\ell=-\frac{1}{\ln (1-p)}$. From this, we can conclude that there exists some \( \ell_0 \) in the interval \( \left(-\frac{1}{\ln(1-p)}, \ell^*\right) \) such that \( g \) is decreasing on \( [2, \ell_0] \) and increasing on \( [\ell_0, \ell^*] \). Therefore, for any $\ell \in [2,\ell^*]$, it holds that $g(\ell) \leq \displaystyle \max \{g(2), g(\ell^*) \}$. 

In order to prove that the right-hand side in~\eqref{bound_1} approaches 0,  we will use the fact that for any $\ell \in [9,\ell^*]$, $g(\ell) \leq \displaystyle \max \{g(9), g(\ell^*) \}$, and also that for any $i\in \{3,\ldots,8 \}$, $e^{g(i)}<e^{g(2)}$. In this way, it remains to show that 
\begin{itemize}
\item $ke^{g(\ell^*)}=o(1)$;
\item $ke^{g(9)}=o(1)$;
\item $e^{g(2)}=o(1)$.
\end{itemize} 
We will show that all these three equalities hold for $p \geq n^{-\theta}$, $\theta \leq 2/7-\varepsilon$ (again, this assumption is weaker that the requirement on $p$ from the statement of Theorem~\ref{th:new}).

First, we obtain that
\begin{align*}
    ke^{g(\ell^*)}&= k \cdot k \left( \frac{k^2e (1-p)}{\ell^* n p} (1-p)^{-\frac{\ell^*}{2}}\right)^{\ell^*}  = k^2\left( \frac{k^2e (1-p)}{\ell^* n p} \left(\frac{1}{1-p} \right)^{\log_{1/(1-p)}\left(\frac{np}{4ek} \right)}\right)^{\ell^*}\\
    &=
    k^2\left( \frac{k (1-p)}{4\ell^*}\right)^{\ell^*}=
    k^2\left( \frac{(2\ln(enp)+o(1))(1-p)}{4(2\ln(np)-2\ln(4ek))}\right)^{\ell^*}
\end{align*}

\begin{align*}
    &=k^2\left( \frac{(2\ln(enp)+o(1))(1-p)}{8\ln(np/k)+O(1)}\right)^{\ell^*}=
    k^2\left(\frac{\ln(np)}{4\ln \left(np^2 \right)}+o(1)\right)^{\ell^*}\\
    & \leq k^2\left(\frac{(1-\theta)}{4(1-2\theta)}+o(1)\right)^{\ell^*}<
    k^2\left(\frac{1}{2}+o(1)\right)^{\ell^*}
    =o(1).
\end{align*}

Next, we get
\begin{equation*}
    ke^{g(9)}= k \cdot k  \left( \frac{k^2e (1-p)}{9 n p} (1-p)^{-\frac{9}{2}}\right)^{9}= O\left(\frac{k^{20}}{n^9 p^9} \right)=O\left(\frac{(\ln n)^{20}}{n^9 p^{29}} \right) = o(1).
\end{equation*}
Lastly, we have
\begin{equation*}
    e^{g(2)}= k \left( \frac{k^2e (1-p)}{2 n p} (1-p)^{-1}\right)^{2}= O\left(\frac{k^{5}}{n^2 p^2} \right)=O\left(\frac{(\ln n)^{5}}{n^2 p^{7}} \right) = o(1).
\end{equation*}

From (\ref{bound_1}), we finally get that for $p \geq n^{-2/7+\varepsilon}$, 
\begin{equation}\label{eq_final_1}
\begin{split}
     \displaystyle \sum_{\ell =2}^{\ell^*} \frac{F_{\ell}}{({\sf E} X_k)^2} & \leq  \displaystyle \sum_{\ell =2}^{\ell^*}  e^{g(\ell)}= \displaystyle \sum_{\ell =2}^{8}  e^{g(\ell)}+\displaystyle \sum_{\ell =9}^{\ell^*}  e^{g(\ell)} \leq 7 e^{g(2)}+k \cdot \max\left\{e^{g(9)}, e^{g(\ell^*)}\right\}=o(1).
\end{split}  
\end{equation}

\subsection{Part 2: $\ell \in \left(\ell^*, \left\lfloor k-\frac{w}{p}\right\rfloor \right]$, $ 1 \ll w= o\left(\sqrt{\ln n}\right)$}
\label{sc:part2}
Here we bound $N(k, \ell, r)\leq k^{k-2}f(k,\ell,r)$, where $f(k,l,r)$, initially defined in~\cite{12}, is an upper bound on the number of trees on $[k]$ that induce $F$ on $[\ell]$ over all choices of a forest $F$ with $r$ edges on $[\ell]$.\\

It is proved in~\cite{12} that $f(k, \ell, r) \left( \frac{1-p}{p}\right)^r \leq (k-\ell)^{k-2}(\ell+1)^{k-\ell-1}$ for all $r\in\{0,\ldots,\ell-1\}$, all $\ell\leq k-\frac{2(1-p)}{p}$, and all $p=p(n)\in(0,1)$. We will not write here the explicit value of $f(k,\ell,r)$ since we will only use this bound on $f(k, \ell, r) \left( \frac{1-p}{p}\right)^r$.
Hence, we obtain for every $\ell \in \left(\ell^*, \left\lfloor k-\frac{w}{p}\right\rfloor\right]$ the following bound for $F_{\ell}$:
\begin{equation*}
    \begin{split}
         F_{\ell} & \leq \binom{n}{k} \binom{k}{\ell} \binom{n-k}{k-\ell} p^{2(k-1)} (1-p)^{2 \binom{k}{2} - \binom{\ell}{2} -2(k-1)} \ell k^{k-2} (k-\ell)^{k-2}(\ell+1)^{k-\ell-1}.
    \end{split}
\end{equation*}
We thus get
\begin{equation}
\begin{split}
    \displaystyle \sum_{\ell =\lfloor\ell^*\rfloor +1}^{\lfloor k-w/p\rfloor} \frac{F_{\ell}}{({\sf E} X_k)^2} & \leq \displaystyle \sum_{\ell =\lfloor\ell^*\rfloor +1}^{\lfloor k-w/p\rfloor} \frac{\binom{n}{k} \binom{k}{\ell} \binom{n-k}{k-\ell} p^{2(k-1)} (1-p)^{2 \binom{k}{2} - \binom{\ell}{2} -2(k-1)} \ell k^{k-2} \displaystyle (k-\ell)^{k-2}(\ell+1)^{k-\ell-1} 
    }{ \binom{n}{k}^2  k^{2(k-2)} p^{2(k-1)}(1-p)^{2\binom{k}{2}-2k+2}} \\ 
    & <  \displaystyle \sum_{\ell =\lfloor\ell^*\rfloor +1}^{\lfloor k-w/p\rfloor}  \frac{\binom{k}{\ell} \binom{n-k}{k-\ell} (1-p)^{-\binom{\ell}{2}}}{\binom{n}{k} k^{k-3}} \displaystyle (k-\ell)^{k-2}(\ell+1)^{k-\ell-1}=
    \displaystyle \sum_{\ell =\lfloor\ell^*\rfloor +1}^{\lfloor k-w/p\rfloor}\hat F_{\ell},
\end{split}
\label{bound_2}
\end{equation}    
where 
$$
\hat{F}_{\ell}=\frac{\binom{k}{\ell} \binom{n-k}{k-\ell} (1-p)^{-\binom{\ell}{2}} (k-\ell)^{k-2} (\ell+1)^{k-\ell-1}}{\binom{n}{k} k^{k-3}}.
$$
The next step is to show that $\hat{F}_{\ell}$ increases in $\ell \in \left(\ell^*,k-\frac{w}{p} \right]$. For this, it is sufficient to prove that $\frac{\hat{F}_{\ell+1}}{\hat{F}_{\ell}}>1$ for all $\ell \in \left(\ell^*,k-\frac{w}{p} \right]$. First, we notice that for any $\ell$ in this interval, the following bounds hold:
\begin{itemize}
\item $(k-\ell)^2\geq \frac{w^2}{p^2}$;
\item $(1-p)^{-\ell} \geq(1-p)^{-\ell^*}= \left( \frac{1}{1-p}\right)^{\frac{\ln \left(\frac{np}{4ek} \right)^2}{\ln(1/[1-p])}}=\frac{n^2p^2}{16 e^2k^2}$;
\item $\left(\frac{k-\ell-1}{k-\ell} \right)^{k-2} =e^{(k-2)\ln(1-\frac{1}{k-\ell})}\geq e^{(k-2)\ln(1-\frac{p}{w})}$.
\end{itemize} 

Therefore, for any $p \geq n^{-1/4+\varepsilon}$, $\varepsilon>0$ (again, this assumption is weaker that the requirement on $p$ from the statement of Theorem~\ref{th:new}), we get
\begin{equation*}
    \begin{split}
         \frac{\hat{F}_{\ell+1}}{\hat{F}_{\ell}} & =\frac{(k-\ell)^2(1-p)^{-\ell}}{(\ell+1)^2(n-2k+\ell+1)}\left(\frac{k-\ell-1}{k-\ell}\right)^{k-2}\left(\frac{\ell+2}{\ell+1} \right)^{k-\ell-2}\\
         & > \frac{1}{(\ell+1)^2}\frac{1}{n} \left(\frac{w}{p} \right)^2 \frac{n^2p^2}{16 e^2k^2} e^{(k-2)\ln(1-\frac{p}{w})}> \frac{n}{16 e^2}\frac{w^2}{k^4}e^{(k-2)\ln(1-\frac{p}{w})} \\
         & > \frac{n}{16 e^2}\frac{w^2 p^4}{3^4 (\ln(npe))^4}e^{(k-2)\ln(1-\frac{p}{w})}\geq \frac{n^{\varepsilon}}{16 e^2}\frac{w^2}{3^4 (\ln(npe))^4}e^{(k-2)\ln(1-\frac{p}{w})},
    \end{split} 
\end{equation*}
and therefore
\begin{equation*}
    \begin{split}
         \frac{\hat{F}_{\ell+1}}{\hat{F}_{\ell}} & > \Omega \left(\frac{w^2}{(\ln n)^4} n^{\varepsilon-O(1/w)}\right)>1.
    \end{split} 
\end{equation*}

In this way, since $\hat{F}_{\ell}$ increases in $\left(\ell^*,k-\frac{w}{p} \right]$, in order to prove that the right-hand side in~\eqref{bound_2} approaches 0, it is sufficient to show that at $\ell=k-\frac{w}{p}$, it holds that $k \cdot \hat{F}_{\ell}$ is $o(1)$. We will prove that this holds for $n^{-1/6+\varepsilon}\leq p=O\left(\frac{1}{\ln n}\right)$, for any $\varepsilon>0$.

Once again, applying Stirling's approximation, we get the following asypmtotics:
\begin{equation*}
    \frac{\binom{k}{\ell} \binom{n-k}{k-\ell}}{\binom{n}{k}} \sim \frac{1}{\sqrt{2\pi}}\frac{k^{2k+1}(n-k)^{2n-2k+1}}{\ell^{1/2+\ell}(k-\ell)^{2k-2\ell+1}n^{n+1/2}(n-2k+\ell)^{n-2k+\ell+1/2}}.
\end{equation*}
Therefore, 
\begin{equation*}
    \begin{split}
            k\cdot \hat{F}_{\ell}& =k \frac{\binom{k}{\ell} \binom{n-k}{k-\ell}}{\binom{n}{k}}\frac{ (1-p)^{-\binom{\ell}{2}} (k-\ell)^{k-2} (\ell+1)^{k-\ell-1}}{ k^{k-3}}\\
           & \sim k \frac{1}{\sqrt{2\pi}}\frac{k^{2k+1}(n-k)^{2n-2k+1}(1-p)^{-\binom{\ell}{2}} (k-\ell)^{k-2} (\ell+1)^{k-\ell-1}}{\ell^{1/2+\ell}(k-\ell)^{2k-2\ell+1}n^{n+1/2}(n-2k+\ell)^{n-2k+\ell+1/2}k^{k-3}}\\
           & = \frac{k^{k+5}(k-\ell)^{3k-2\ell-3}(n-k)^{2n-2k+1}}{(1-p)^{\binom{\ell}{2}}\ell^{\ell+1/2}}n^{-n-1/2}(n-2k+\ell)^{2k-\ell-n-1/2}(\ell+1)^{k-\ell-1}\\
           & = e^{h(\ell)},
    \end{split}
\end{equation*}
where 
\begin{equation*}
    \begin{split}
        h(\ell)=&(k+5)\ln k+(3k-2\ell-3)\ln (k-\ell)+(2n-2k+1)\ln(n-k)+(-\ell-1/2)\ln \ell\\
    &+(-n-1/2)\ln n+(2k-\ell-n-1/2)\ln(n-2k+\ell)-\binom{\ell}{2}\ln(1-p)\\
    &+(k-\ell-1)\ln(\ell+1).
    \end{split}
\end{equation*}
Letting $\ell=k-\frac{w}{p}$, we get
\begin{equation*}
    \begin{split}
        h(\ell)=&(k+5)\ln k+(2\frac{w}{p}+k-3)\ln (\frac{w}{p})+(2n-2k+1)\ln(n-k)+(\frac{w}{p}-k-1/2)\ln (k-\frac{w}{p})\\
    &+(-n-1/2)\ln n+(\frac{w}{p}+k-n-1/2)\ln(n-\frac{w}{p}-k)+(\frac{w}{p}-1)\ln(k-\frac{w}{p}+1)\\
    &-\frac{k^2}{2}\ln(1-p)+k\frac{w}{p} \ln(1-p)-\frac{w^2}{2p^2}\ln(1-p) +\frac{k}{2}\ln(1-p)-\frac{w}{2p}\ln(1-p)\\
    =& k \left[\ln k+\ln(\frac{w}{p})-2\ln n-\ln k+\ln n-\frac{k}{2}\ln(1-p)+\frac{w}{p}\ln(1-p) \right]\\
    &+\frac{w}{p}\left[2\ln(\frac{w}{p})+2\ln k+\ln n -\frac{w}{2p}\ln(1-p) \right]+O(\ln k)+O(\ln(\frac{w}{p}))+O(\ln n)\\
     =& k \left[\ln(\frac{w}{p})-\ln n+\frac{kp}{2}+\frac{kp^2}{4}+O(kp^3)-w+O(wp) \right]\\
    &+\frac{w}{p}\left[2\ln(\frac{w}{p})+2\ln k+\ln n +\frac{w}{2}+O(wp)\right]+O(\ln k)+O(\ln(\frac{w}{p}))+O(\ln n).
    \end{split}
\end{equation*}
Substituting the expression from~\eqref{value_k} for $k$ in some of the summands, we obtain
\begin{equation*}
    \begin{split}
    h(\ell)=& k \left[\ln(\frac{wnp}{pn})+p\ln(np)+O(kp^3)-w+O(wp) \right]\\
    &+\frac{w}{p}\left[-4\ln p +2\ln w+2\ln \ln (np)+\ln n +\frac{w}{2}+O(wp)\right]+O(\ln k)+O(\ln(\frac{w}{p}))+O(\ln n)\\
    =& k \left[\ln w+p\ln (np)-w+O(kp^3)+O(wp) \right]\\
    &+\frac{w}{p}\left[\ln(\frac{n}{p^4})+2\ln w+\frac{w}{2}+2\ln \ln (np)+O(wp)\right]+O(\ln k)+O(\ln(\frac{w}{p}))+O(\ln n)\\
    =& k \left[\ln w+p\ln (np)-w\right]+\frac{w}{p}\ln\left(\frac{n}{p^4}\right)+O(p(\ln n)^2)+O(w \ln n)+2\frac{w}{p}\ln w+\frac{w^2}{2p}\\
    &+2\frac{w}{p}\ln \ln (np)+O(w^2) +O(\ln k)+O(\ln(\frac{w}{p}))+O(\ln n)\\
    =& \underbrace{k \ln w}_{(\romannumeral 1)}+\underbrace{kp\ln (np)}_{(\romannumeral 2)}-\underbrace{kw}_{(\romannumeral 3)}+\underbrace{\frac{w}{p}\ln\left(\frac{n}{p^4}\right)}_{(\romannumeral 4)}+O(p(\ln n)^2)+O(w \ln n)+O(\frac{w^2}{p})+\frac{2w}{p}\ln \ln (np).
    \end{split}
\end{equation*}
Observe that if \( p =\Omega \left( \frac{1}{\sqrt{\ln n}}\right) \), then the dominant asymptotics in \( h(\ell) \) will be determined by \( (\romannumeral 2) \). Since \( (\romannumeral 2) \to \infty \), it follows that \( k \cdot \hat{F}_{\ell} \) at \( \ell = k - \frac{w}{p} \) is not \( o(1) \). Nevertheless, for $n^{-\theta}\leq p=O\left(\frac{1}{\ln n}\right)$, for some $\theta>0$, the dominant asymptotics in \( h(\ell) \) will be determined by \( (\romannumeral 3) \) and \( (\romannumeral 4) \), for which, as soon as $p\geq n^{-1/6+\varepsilon}$, we have
\begin{equation*}
    -kw+\frac{w}{p}\ln\left(\frac{n}{p^4}\right)=-2\frac{w}{p}\ln (np)+\frac{w}{p}\ln\left(\frac{n}{p^4}\right)+O\left(\frac{w}{p}\right)=-\frac{w}{p}\ln \left(np^6 \right)+O\left(\frac{w}{p}\right)\to -\infty.
\end{equation*}

From (\ref{bound_2}) and according to the above remarks, we finally get that for $n^{-1/6+\varepsilon}\leq p=O\left(\frac{1}{\ln n}\right)$, $\varepsilon>0$, it holds that

\begin{equation}\label{eq_final_2}
    \displaystyle \sum_{\ell =\lfloor\ell^*\rfloor +1}^{\lfloor k-w/p\rfloor} \frac{F_{\ell}}{({\sf E} X_k)^2}  \leq \displaystyle \sum_{\ell =\lfloor\ell^*\rfloor +1}^{\lfloor k-w/p\rfloor}\hat F_{\ell} \leq  k \cdot \hat{F}_{\ell}|_{\ell=k-\frac{w}{p}}=o(1).
\end{equation}

\subsection{Part 3: $\ell \in \left(\left\lfloor k-\frac{w}{p}\right\rfloor, k- \frac{1}{2p}\right]$, $1 \ll w= o\left(\sqrt{\ln n}\right)$}
\label{sc:part3}
Here we bound $N(k, \ell, r)\leq \varphi(\ell,r)[f(k,\ell,r)]^2$, where $\varphi(\ell,r)$ is the number of forests $F$ on $[\ell]$ with $r$ edges, and $f(k,l,r)$ is an upper bound on the number of trees on $[k]$ that induce $F$ on $[\ell]$ over all choices of $F$. \\

From \cite{15} the number of rooted forests in $[n]$ with $m$ trees is given by $\binom{n}{m}mn^{n-m+1}$. Then, since a forest on $[\ell]$ with $r$ edges has $\ell-r$ trees, we get $\varphi(\ell,r) \leq \binom{\ell}{\ell-r} (\ell -r)\ell^{r-1}$. Thus we obtain for every $\ell \in \left(\left\lfloor k-\frac{w}{p}\right\rfloor, k- \frac{1}{2p}\right]$ the following bound for $F_{\ell}$:
\begin{equation}
    \begin{split}
         F_{\ell} & \leq \binom{n}{k} \binom{k}{\ell} \binom{n-k}{k-\ell} p^{2(k-1)} (1-p)^{2 \binom{k}{2} - \binom{\ell}{2} -2(k-1)} \ell \max_{r \in\{0,\ldots,\ell-1\}} H(k,\ell,r),
    \end{split}
    \label{fl_part3}
\end{equation}
where $H(k, \ell, r) = \left(\frac{1-p}{p} \right)^r \binom{\ell}{\ell-r} (\ell-r)\ell^{r-1}[f(k, \ell, r)]^2$.

Our next step is to find for which values of $r$, the function $H(k, \ell, r)$ is maximal. For this, we will use the values of $f(k, \ell, r)$ established in \cite{12} for different intervals of $r$. In Section \ref{para_1} we study $H(k, \ell, r)$ for $r \in \left[\ell(1-1/e), \ell-1\right)$, in Section \ref{para_2} --- for $r \in \left[\ell/2, \ell(1-1/e)\right)$, and in Section \ref{para_3} --- for $r \in \left[0, \ell/2\right)$. In these sections, we will show that $H(k, \ell, r)$ increases in $\left[0, \ell(1-1/e)\right)$, then in $r=\ell(1-1/e)$ it has an increasing jump, and finally, for some $r^*$, it increases in $\left[\ell(1-1/e), r^*\right]$ and decreases in $\left[r^*, \ell-1\right)$. In this way, we will show that $H(k, \ell, r)$ is maximal in $r=r^*$. In Section~\ref{para_4}, we find the value of $r^*$.

\subsubsection{$H(k, \ell, r)$ in $r \in \left[\ell(1-1/e), \ell-1\right)$.}
\label{para_1}
From \cite{12} we have 
$$
f(k,\ell,r)=\left(\frac{\ell}{\ell-r} \right)^{\ell-r}(\ell+1)^{k-\ell-1}(k-\ell)^{k-r-2}
$$ 
for $r \in \left[\ell(1-1/e), \ell-1\right)$. We thus get
\begin{equation}\label{right_lim_fin}
    H(k,\ell,r)=\left(\frac{1-p}{p} \right)^r \binom{\ell}{\ell-r} (\ell-r)\ell^{r-1} \left(\frac{\ell}{\ell-r} \right)^{2(\ell-r)}(\ell+1)^{2(k-\ell-1)}(k-\ell)^{2(k-r-2)}.
\end{equation}
It follows that
\begin{equation}
    \begin{split}
        \frac{H(k,\ell,r+1)}{H(k,\ell,r)}=\frac{(1-p)}{p(k-\ell)^2\ell} \cdot \frac{(\ell-r)^{2(\ell-r)}}{(\ell-r-1)^{2(\ell-r)-3} (r+1)}=\frac{(1-p)}{p(k-\ell)^2\ell} \cdot \varphi(r),
    \end{split}
    \label{eq_H}
\end{equation}
where $\varphi(r)=\frac{(\ell-r)^{2(\ell-r)}}{(\ell-r-1)^{2(\ell-r)-3} (r+1)}$. Notice that $\varphi(\ell(1-1/e))=\frac{1}{e^2}\left(\frac{\ell}{\ell-e}\right)^{2\ell/e}\frac{(\ell-e)^3}{e\ell-\ell+e}\sim \ell^2$ and $\varphi(\ell-2)=\frac{16}{\ell-1}$. Therefore, since $w= o\left(\sqrt{\ln n}\right)$, we get
\begin{equation}
    \begin{split}
        \frac{H(k,\ell,r+1)}{H(k,\ell,r)}\bigg|_{r=\ell(1-1/e)}\sim \frac{(1-p)}{p(k-\ell)^2\ell}\ell^2=\Theta \left(\frac{\ell}{p(k-\ell)^2} \right)=\Omega \left(\frac{\ell p}{w^2} \right)=\Omega \left(\frac{\ln (np)}{w^2} \right) \to \infty,
    \end{split}
    \label{h_1}
\end{equation}
and
\begin{equation}
    \begin{split}
        \frac{H(k,\ell,r+1)}{H(k,\ell,r)}\bigg|_{r=\ell-2}=\frac{(1-p)}{p(k-\ell)^2\ell}\frac{16}{(\ell-1)}=O\left(\frac{1}{p\left(\frac{1}{2p} \right)^2 \ell^2}\right)=O\left(\frac{p}{\ell^2}\right) \to 0.
    \end{split}
    \label{h_2}
\end{equation}
Moreover, since $\ln\varphi(r)=2(\ell-r)\ln (\ell-r)-(2(\ell-r)-3)\ln(\ell-r-1)-\ln(r+1)$, then for all $r \in \left[\ell(1-1/e), \ell-1\right)$ we get
\begin{equation}
    \begin{split}
       \frac{\partial}{\partial r}\ln\varphi(r) & =-\frac{2\ell}{\ell-r}+\frac{2r}{\ell-r}-2\ln(\ell-r)+\frac{2\ell}{\ell-r-1}+2\ln(\ell-r-1)-\frac{2r+3}{\ell-r-1}-\frac{1}{r+1}\\
       & = -\frac{2(\ell-r)}{\ell-r}+\frac{2(\ell-r-1)-1}{\ell-r-1}+2\ln\left(\frac{\ell-r-1}{\ell-r} \right)-\frac{1}{r+1}\\
       & =-\frac{1}{\ell-r-1}+2\ln\left(\frac{\ell-r-1}{\ell-r} \right)-\frac{1}{r+1}<0.
    \end{split}
    \label{dev_H}
\end{equation}
From (\ref{dev_H}) it follows that $\frac{H(k,\ell,r+1)}{H(k,\ell,r)}$ strictly decreases in $r \in \left[\ell(1-1/e), \ell-1\right)$. This, together with (\ref{h_1}) and (\ref{h_2}), gives that there exists an unique point $r^*$, being the solution of $\frac{H(k,\ell,r^*+1)}{H(k,\ell,r^*)}=1$, such that $H(k,\ell,r)$ increases in $r \in \left[\ell(1-1/e), r^*\right]$ and decreases in $r \in \left[r^*, \ell-1\right)$.

\subsubsection{$H(k, \ell, r)$ in $r \in [\ell/2, \ell(1-1/e))$.}
\label{para_2}

From \cite{12} we have $f(k,\ell,r)=3^{2r-\ell} 2^{2\ell-3r} (\ell+1)^{k-\ell-1}(k-\ell)^{k-r-2}$ for $r \in [\ell/2, \ell(1-1/e))$. We thus get
\begin{equation}\label{right_lim}
    H(k,\ell,r)=\left(\frac{1-p}{p} \right)^r \binom{\ell}{\ell-r} (\ell-r)\ell^{r-1} 3^{4r-2\ell} 2^{4\ell-6r} (\ell+1)^{2(k-\ell-1)}(k-\ell)^{2(k-r-2)}.
\end{equation}

In this way, for all $r \in [\ell/2, \ell(1-1/e))$, it holds 
\begin{equation*}
    \begin{split}
        \frac{H(k,\ell,r+1)}{H(k,\ell,r)}& =\frac{81(1-p)\ell}{64p(k-\ell)^2}\frac{\ell-r-1}{r+1}\geq \frac{(1-p)}{p(k-\ell)^2}(\ell-r-1)\geq \frac{1-p}{p}\frac{p^2}{w^2}\left(\frac{\ell}{e}-1\right)\\
        &=\frac{(1-p)p}{w^2}\left(\frac{\ell}{e}-1\right) =\Theta\left(\frac{p \ell}{w^2} \right)=\Theta\left(\frac{\ln n}{w^2
        } \right) \to \infty.
    \end{split}
    \label{eq_H_2}
\end{equation*}

Therefore $H(k, \ell, r)$ increases in $r \in [\ell/2, \ell(1-1/e))$.

\subsubsection{$H(k, \ell, r)$ in $r \in \left[0, \ell/2\right)$.}
\label{para_3} 
From \cite{12} we have $f(k,\ell,r)= 2^{r} (\ell+1)^{k-\ell-1}(k-\ell)^{k-r-2}$ for $r \in \left[0, \ell/2\right)$. We thus get

\begin{equation}\label{left_lim}
    H(k,\ell,r)=\left(\frac{1-p}{p} \right)^r \binom{\ell}{\ell-r} (\ell-r)\ell^{r-1} 2^{2r} (\ell+1)^{2(k-\ell-1)}(k-\ell)^{2(k-r-2)}.
\end{equation}

In this way, for all $r \in \left[0, \ell/2\right)$, it holds 
\begin{equation*}
    \begin{split}
        \frac{H(k,\ell,r+1)}{H(k,\ell,r)}& = \frac{(1-p)4\ell}{p (k-\ell)^2} \frac{(\ell-r-1)}{r+1} \geq \frac{(1-p)}{p}\frac{4}{(k-\ell)^2}\left(\frac{\ell}{2}-1 \right)\geq  \frac{(1-p)}{p} \frac{p^2}{w^2}(\ell-2)\\
        & =\frac{p(1-p)(\ell-2)}{w^2}=\Theta \left(\frac{p\ell}{w^2} \right)=\Theta\left(\frac{\ln n}{w^2} \right) \to \infty.
    \end{split}
    \label{eq_H_2}
\end{equation*}

Therefore $H(k, \ell, r)$ increases in $r \in \left[0, \ell/2\right)$.\\

The only point remaining concerns the behaviour of $H(k, \ell, r)$ at the end-points of the considered intervals, i.e. at $r\in \{\ell/2,\ell(1-1/e) \}$. It follows easily that the the left and right-hand limits of $H(k, \ell, r)$ at $r=\ell/2$ coincide. The left-hand limit of $H(k, \ell, r)$ at this points is given by (\ref{left_lim}) and the right-hand limit by (\ref{right_lim}). In both cases it equals
\begin{equation}\label{equalities}
    \left(\frac{1-p}{p} \right)^{\ell/2} \binom{\ell}{\ell/2} \frac{\ell}{2}\ell^{\frac{\ell}{2}-1} 2^{\ell} (\ell+1)^{2(k-\ell-1)}(k-\ell)^{2(k-2-\ell/2)}.
\end{equation}

On the other hand, (\ref{right_lim}) gives the following left-hand limit of $H(k, \ell, r)$ at $r=\ell(1-1/e)$

\begin{equation}\label{left_1}
    \left(\frac{1-p}{p} \right)^{\ell(1-1/e)} \binom{\ell}{\ell/e} \frac{\ell}{e}\ell^{\ell(1-1/e)-1}(\ell+1)^{2(k-\ell-1)}(k-\ell)^{2(k-\ell(1-1/e)-2)} 3^{2\ell(1-\frac{2}{e})} 2^{-2\ell(\frac{3}{e}+1)},
\end{equation}

and (\ref{right_lim_fin}) gives the following right-hand limit

\begin{equation}\label{right_1}
    \left(\frac{1-p}{p} \right)^{\ell(1-1/e)} \binom{\ell}{\ell/e} \frac{\ell}{e}\ell^{\ell(1-1/e)-1}(\ell+1)^{2(k-\ell-1)}(k-\ell)^{2(k-\ell(1-1/e)-2)} e^{\frac{2\ell}{e}}.
\end{equation}

We can see that (\ref{left_1}) is less than (\ref{right_1}). We thus get that $H(k, \ell, r)$ has an increasing jump at $r=\ell(1-1/e)$, and we can finally conclude that $H(k, \ell, r)$ is maximal at $r=r^*$. 

\subsubsection{Finding $r^*$}
\label{para_4}

The next step is to compute the value of $r^*$. For this, we find $r$, such that (\ref{eq_H}) equals $1+o(1)$. Let us recall that here we consider $\ell \in \left(\left\lfloor k-\frac{w}{p}\right\rfloor, k- \frac{1}{2p}\right]$. For simplicity of notation, we will write $\ell=k-\frac{\beta}{p}$, where $\beta \in \left[\frac{1}{2},w \right]$. Let us find $r^*$ in the form $r^*=\ell-\frac{\lambda}{p}$. We thus get

\begin{equation*}
    \begin{split}
        \frac{H(k,\ell,r+1)}{H(k,\ell,r)} &=\frac{(1-p)}{p(k-\ell)^2\ell} \frac{(\ell-r)^{2(\ell-r)}}{(\ell-r-1)^{2(\ell-r)-3} (r+1)}\\
        & =e^{\ln(1-p)-\ln p-2 \ln (k-\ell) - \ln \ell -\ln (r+1) +2(\ell-r)\ln(\ell-r)-[2(\ell-r)-3]\ln(\ell-r-1)}\\
        & =e^{-p+O(p^2)-\ln p -2 \ln \left(\frac{\beta}{p}\right)-\ln \ell -\ln(r+1)-2(\ell-r)\ln\left(1-\frac{1}{\ell-r} \right)+3\ln(\ell-r-1)}\\
        & = e^{-2\ln \beta+\ln p -\ln \ell -\ln(\ell-\frac{\lambda}{p}+1)-2(\ell-(\ell-\frac{\lambda}{p}))\ln\left(1-\frac{1}{\ell-(\ell-\frac{\lambda}{p})} \right)+3\ln(\ell-(\ell-\frac{\lambda}{p})-1)+o(1)}\\
        & =e^{-2\ln \beta+\ln p -\ln \ell-\ln\left(\ell-\frac{\lambda}{p} \right)-2\frac{\lambda}{p}\ln \left( 1-\frac{p}{\lambda}\right)+3\ln\left(\frac{\lambda}{p}-1 \right)+o(1)}\\
         & =e^{-2\ln \beta+\ln p -2 \ln \ell+2+3\ln\left(\frac{\lambda}{p} \right)+o(1)}\\
         & = e^{\ln \left(\lambda^3\left(\frac{e}{\beta \ell p}\right)^2 \right)+o(1)} \to 1, {\text{ if }} \lambda=\left(\frac{\beta \ell p}{e}\right)^{2/3}.
    \end{split}
    \label{eq_H_tozero}
\end{equation*}

Therefore $r^*=\ell-\frac{1}{p}\left(\frac{\beta \ell p}{e}\right)^{2/3}$. Notice that $\lambda=\Theta(\beta^{2/3}(\ln n)^{2/3})=o(\ln n)$.\\

\subsubsection{Completing the proof}

Having finished the analysis of $H(k, \ell, r)$, and since $r^* \in \left[\ell(1-1/e), \ell-1\right)$, we can use (\ref{right_lim_fin}) to obtain the following bound on $H(k, \ell, r)$. Here we also use the known bound $\binom{n}{k}\leq \left(\frac{en}{k}\right)^k$.
\begin{equation}
    \begin{split}
        &\max_{r} H(k,\ell,r)  \leq H(k,\ell,r^*)= \\
        & = \left(\frac{1-p}{p} \right)^{r^*} \binom{\ell}{\ell-r^*} (\ell-r^*)\ell^{r^*-1} \left(\frac{\ell}{\ell-r^*} \right)^{2(\ell-r^*)}(\ell+1)^{2(k-\ell-1)}(k-\ell)^{2(k-r^*-2)}\\
        & \leq \left(\frac{1-p}{p} \right)^{r^*} \frac{e^{\ell-r^* }\ell^{\ell-r^*}}{(\ell-r^*)^{\ell-r^*}}(\ell-r^*)\ell^{r^*-1}\frac{\ell^{2(\ell-r^*)}}{(\ell-r^*)^{2(\ell-r^*)}}(\ell+1)^{2(k-\ell-1)}(k-\ell)^{2(k-r^*-2)} \\
        & = \left(\frac{1-p}{p} \right)^{r^*} e^{\ell-r^*}\ell^{3\ell-2r^*-1}(\ell-r^*)^{3r^*-3\ell+1}(\ell+1)^{2(k-\ell-1)}(k-\ell)^{2(k-r^*-2)}.
    \end{split}
    \label{upp_boundd}
\end{equation}
Before moving to the final step, let us modify the bound in (\ref{bound_11}) in the following way. We get that for some constant $c<1$ and for all sufficiently large $n$, it holds
\begin{equation}
    \begin{split}
      \frac{\binom{k}{\ell} \binom{n-k}{k-\ell}}{\binom{n}{k}} & < c\left(\frac{k^2}{\ell n} \right)^{\ell}
\left(\frac{k}{k-\ell}\right)^{2(k-\ell)}
\left(\frac{n}{n-k+\ell}\right)^{\ell}
\left(\frac{n-k}{n-k+\ell}\right)^{n-k}\\
&\leq c\left(\frac{k^2}{\ell n} \right)^{\ell} \left(\frac{k}{k-\ell}\right)^{2(k-\ell)}
\exp\left[\frac{\ell(k-\ell)}{n-k+\ell}-\frac{\ell(n-k)}{n-k+\ell}\right]\\
& = c\left(\frac{k^2}{\ell n} \right)^{\ell} \left(\frac{k}{k-\ell}\right)^{2(k-\ell)} \exp\left[-\ell+\frac{\ell k}{n-k+\ell}\right]\\
& \leq c \left( \frac{k^2}{\ell ne} \right)^{\ell}\left(\frac{k}{k-\ell}\right)^{2(k-\ell)}. 
    \end{split}
    \label{bound_12}
\end{equation}
Thus, by applying~\eqref{fl_part3} and~\eqref{upp_boundd}, and denoting
$$
H(\ell):=\frac{\binom{k}{\ell} \binom{n-k}{k-\ell} \ell\cdot\left(\frac{1-p}{p}\right)^{r^*}}{\binom{n}{k} k^{2(k-2)}(1-p)^{\binom{\ell}{2}}} e^{\ell-r^*}\ell^{3\ell-2r^*-1}(\ell-r^*)^{3(r^*-\ell)+1}(\ell+1)^{2(k-\ell-1)}(k-\ell)^{2(k-r^*-2)},
$$
we obtain that
\begin{equation}
    \begin{split}
        & \displaystyle \sum_{\ell =\lfloor k-w/p\rfloor}^{k-\frac{1}{2p}} \frac{F_{\ell}}{({\sf E} X_k)^2} \leq  \displaystyle \sum_{\ell =\lfloor k-w/p\rfloor}^{k-\frac{1}{2p}} H(\ell).
    \end{split}\label{sum_part_3}
\end{equation}
Our approach here will be to upper bound (\ref{sum_part_3}) by \( k \cdot H(\ell) \). Since we do not have an explicit upper bound for \( H(\ell) \), we will prove that for any \( \ell \in \left(\left\lfloor k-\frac{w}{p}\right\rfloor, k- \frac{1}{2p}\right] \), we have \( k \cdot H(\ell) = o(1) \) for $n^{-1/3+\varepsilon}\leq p<\frac{1}{2\ln n}$ and any $\varepsilon>0$.

Therefore, by (\ref{bound_12}) and (\ref{sum_part_3}) we have that, for some $\ell \in \left(\left\lfloor k-\frac{w}{p}\right\rfloor, k- \frac{1}{2p}\right]$,
\begin{equation}
    \begin{split}
        & \displaystyle \sum_{\ell =\lfloor k-w/p\rfloor}^{k-\frac{1}{2p}} \frac{F_{\ell}}{({\sf E} X_k)^2} \leq \displaystyle \sum_{\ell =\lfloor k-w/p\rfloor}^{k-\frac{1}{2p}} H(\ell) \stackrel{\text{fix }\ell}\leq \\
    & \leq \frac{\binom{k}{\ell} \binom{n-k}{k-\ell} (1-p)^{-\binom{\ell}{2}} }{\binom{n}{k} k^{2(k-2)-1}}\left(\frac{1-p}{p} \right)^{r^*} e^{\ell-r^*}\ell^{3\ell-2r^*}(\ell-r^*)^{3(r^*-\ell)+1}(\ell+1)^{2(k-\ell-1)}(k-\ell)^{2(k-r^*-2)}\\
    & \leq \left( \frac{k^2}{\ell ne} \right)^{\ell}\left(\frac{k}{k-\ell}\right)^{2(k-\ell)}  \frac{(1-p)^{-\binom{\ell}{2}}(\ell+1)^{2(k-\ell-1)}}{k^{2k-5}(k-\ell)^{-2(k-r^*-2)}}\left(\frac{1-p}{p} \right)^{r^*} e^{\ell-r^*}\ell^{3\ell-2r^*}(\ell-r^*)^{3r^*-3\ell+1}\\
    & =k^5 \ell^{2\ell-2r^*}n^{-\ell}(k-\ell)^{2(\ell-r^*-2)}e^{-r^*}\left(\frac{1-p}{p} \right)^{r^*}(\ell-r^*)^{3r^*-3\ell+1}(\ell+1)^{2(k-\ell-1)}(1-p)^{\frac{\ell}{2}-\frac{\ell^2}{2}}\\
    & = e^{g(\ell)},
    \end{split}\label{e^g}
\end{equation}
where
\begin{equation*}
    \begin{split}
         g(\ell) = & 5\ln k +[2(\ell-r^*)]\ln \ell-\ell \ln n +2(\ell-r^*-2)\ln (k-\ell) -r^*+r^* \ln(1-p)-r^* \ln p\\
         & +\left[\frac{\ell}{2}-\frac{\ell^2}{2}\right]\ln (1-p)+[3(r^*-\ell)+1]\ln(\ell-r^*)+2(k-\ell-1)\ln(\ell+1).
    \end{split}
    \label{g_l_1}
\end{equation*}

Our goal is to show $g(\ell) \to -\infty$. In what follows, we write $k-\frac{\beta}{p}$ instead of $\ell$, $k-\frac{\beta}{p}-\frac{\lambda}{p}$ instead of $r^*$, and $g(\beta)$ instead of $g(\ell=k-\beta/p)$. We thus get 
\begin{equation}
    \begin{split}
         g(\beta) = & 5\ln k +\left[\frac{2\lambda}{p}\right]\ln \left( k-\frac{\beta}{p}\right)-\left( k-\frac{\beta}{p}\right) \ln n +\left(\frac{2\lambda}{p}-4\right)\ln \left(\frac{\beta}{p} \right) \\
         &+\left(k-\frac{\beta}{p}-\frac{\lambda}{p} \right)[\ln(1-p)-\ln p -1] +\left[\frac{k}{2}-\frac{\beta}{2p}-\frac{k^2}{2}+\frac{k\beta}{p}-\frac{\beta^2}{2p^2} \right] \ln(1-p)\\
         &+\left[1-\frac{3\lambda}{p} \right]\ln \left( \frac{\lambda}{p}\right) +2\left[\frac{\beta}{p}-1\right]\ln \left(k-\frac{\beta}{p}+1 \right).
    \end{split}
    \label{g_l_1}
\end{equation}

Since 
$$5\ln k-4\ln(\beta/p)+\left[\frac{k}{2}-\frac{\beta}{2p}\right]\ln(1-p)+\ln(\lambda/p)-2\ln\left(k-\frac{\beta}{1}+1\right)=O(\ln n),$$ 
and rearranging the order of the summands in (\ref{g_l_1}) to simplify the asymptotic analysis, we obtain
\begin{equation*}
    \begin{split}
         g(\beta) = & -k \ln n +k[\ln(1-p)-\ln p -1]-\frac{k^2}{2}\ln(1-p)  +\frac{2\lambda}{p}\ln \left( k-\frac{\beta}{p}\right)+\frac{\beta}{p} \ln n  \\
         &+\frac{2\lambda}{p}\ln \left(\frac{\beta}{p} \right)-\left(\frac{\beta}{p}+\frac{\lambda}{p} \right)[\ln(1-p)-\ln p -1]+\left[\frac{k\beta}{p}-\frac{\beta^2}{2p^2} \right] \ln(1-p)\\
         &-\frac{3\lambda}{p}\ln \left( \frac{\lambda}{p}\right)+2\frac{\beta}{p}\ln \left(k-\frac{\beta}{p}+1 \right)+O(\ln n).
\end{split}
\end{equation*}
Therefore,
\begin{equation*}
    \begin{split}
          g(\beta) =&  -k\left[\ln n+\ln p-\frac{kp}{2}+1-\frac{kp^2}{4}+p\right]\\
         &+ \frac{\beta}{p}\left[ \ln n +\ln p+1+2\ln k+k \ln (1-p)-\frac{\beta}{2p}\ln(1-p)-2\frac{\beta}{kp}+O\left(\frac{\beta^2}{(\ln n)^2}\right)\right]\\
         & + \frac{\lambda}{p}\left[2\ln k +2\ln p-\frac{2\beta}{kp}+1-3\ln \lambda+o(1)\right]+O(\ln n)\\
          =&  -k\left[\ln (np)-\ln (npe)+1-p\frac{\ln(npe)}{2}+o(1)\right]\\
         &+ \frac{\beta}{p}\left[ \ln (np) +2\ln k-2\ln(npe)+\frac{\beta}{2}+O\left(\frac{\beta^2}{(\ln n)^2}\right)+O(p\ln n)+O(\beta p)+O(1)\right]\\
         & + \frac{\lambda}{p}\left[-2\ln p +2\ln \ln n +2\ln p-\frac{\beta}{\ln (np)}-3\ln \lambda+O(1)\right]+O(\ln n),
\end{split}
\end{equation*}
where several terms involving $k$ cancel out, leading to the expression
\begin{equation*}
    \begin{split}
         g(\beta) = &  -k\left[-p\frac{\ln(npe)}{2}+o(1)\right]\\
         &+ \frac{\beta}{p}\left[ -2\ln p+2\ln \ln n-\ln(npe)+\frac{\beta}{2}+O\left(\frac{\beta^2}{(\ln n)^2}\right)+O(p\ln n)+O(\beta p)+O(1)\right]\\
         & + \frac{\lambda}{p}\left[2\ln \ln n-\frac{\beta}{\ln (np)}-3\ln \lambda+O(1)\right]+O(\ln n)\\
         =&  -k\left[-p\frac{\ln(npe)}{2}+o(1)\right]+\frac{\beta}{p}\left[ -\ln(np^3)+2\ln \ln n+o(\ln n)+O(1)\right]\\
         & + \frac{\lambda}{p}\left[2\ln \ln n+O(\ln \lambda)\right]+O(\ln n),
    \end{split}
    \label{g_l_2}
\end{equation*}

and this finally leads to
\begin{equation*}
    \begin{split}
         g(\beta) = & \underbrace{kp\frac{\ln(npe)}{2}}_{\romannumeral 1}-\underbrace{\frac{\beta}{p}\ln(np^3)}_{\romannumeral 2}+\underbrace{\frac{\lambda}{p}2\ln \ln n}_{\romannumeral 3}+  \frac{\beta}{p}\left[ 2\ln \ln n+o(\ln n)+O(1)\right]\\
         & +O\left(\frac{\lambda}{p}\ln \lambda\right)+O(\ln n)+o(k).
    \end{split}
    \label{g_l_2}
\end{equation*}

Let us provide a description of the asymptotics of $g(\beta)$ as follows:\\

\textbf{1. } For $p=\frac{c}{\ln n}$, $c>0$, the dominant asymptotics in \( g(\beta) \) will be determined by \( (\romannumeral 1) \) and  \( (\romannumeral 2) \)\footnote{Observe that if \( p = \Omega\left(\frac{1}{\sqrt{\ln n}}\right) \), then the dominant asymptotics in \( g(\beta) \) will be determined by \( (\romannumeral 1) \). Since \( (\romannumeral 1) \to \infty \), it follows that the right-hand side of (\ref{e^g}) is not $o(1)$.}, for which we have
\begin{equation*}
\begin{split}
    kp\frac{\ln(npe)}{2}-\frac{\beta}{p}\ln(np^3)&=\frac{kp}{2}\ln n-\frac{\beta}{p}\ln n+O(\beta \ln n \ln \ln n)\\
    &=\left(1-\frac{\beta}{c} \right)\ln n+O(\beta \ln n \ln \ln n)\to -\infty,
\end{split}
\end{equation*}
for $c<\frac{1}{2}$. Hence, we get that right-hand side of (\ref{e^g}) is $o(1)$ for $p<\frac{1}{2\ln n}$. \\

\textbf{2. } Let $n^{-1/3+\varepsilon}\ll p\ll\ \frac{1}{\ln n}$. Denote $p=n^{-\theta}$, $\theta=\theta(n)>0$. The dominant asymptotics in \( g(\beta) \) will be determined by \( (\romannumeral 2) \) and  \( (\romannumeral 3) \), for which we have
\begin{equation*}
    \begin{split}
        -\frac{\beta}{p}\ln(np^3)+2\frac{\lambda}{p}\ln \ln n& =\frac{1}{p}\left[-\beta (1-3\theta)\ln n+2\left(\frac{\beta \ell p}{e} \right)^{2/3}\ln \ln n \right]\\
        & =\frac{1}{p}\left[-\beta (1-3\theta)\ln n+2\left(\frac{\beta}{e}(kp-\beta) \right)^{2/3}\ln \ln n \right]\\
        & \leq \frac{1}{p}\left[-\beta (1-3\theta)\ln n+2\left(\frac{\beta}{e}(2 \ln(np)) \right)^{2/3}\ln \ln n \right]\\
        & = \frac{1}{p}\left[-\beta (1-3\theta)\ln n+2 \left(\frac{2\beta (1-\theta) }{e}\right)^{2/3} (\ln n)^{2/3}\ln \ln n \right]\\
        & = \beta \frac{\ln n}{p}\left[- (1-3\theta)+2 \left(\frac{2 (1-\theta) }{e}\right)^{2/3} \frac{\ln \ln n}{(\beta \ln n)^{1/3}} \right]\\
        & = \beta \frac{\ln n}{p}\left[- (1-3\theta)+o(1)\right] \to -\infty, \ \text{ for } \theta <\frac{1}{3}. 
    \end{split}
\end{equation*}


We may conclude that, for any $n^{-1/3+\varepsilon} \leq p <\frac{1}{2\ln n}$, $\varepsilon>0$, we have $g(\ell) \to -\infty$, and so
\begin{equation}\label{eq_final_3}
        \displaystyle \sum_{\ell =\lfloor k-w/p\rfloor}^{k-\frac{1}{2p}} \frac{F_{\ell}}{({\sf E} X_k)^2} \leq e^{g(\ell)}=o(1),
\end{equation}
as required.

\subsection{Part 4: $\ell \in \left(k-\frac{1}{2p}, k-1\right]$}
\label{sc:part4}
Here we bound $N(k, \ell, r)\leq k^{k-2}f(k,\ell,r)$, where $f(k,l,r)$ is the upper bound on the number of trees on $[k]$ that induce $F$ on $[\ell]$ over all choices of a forest $F$ with $r$ edges on $[\ell]$, defined in~\cite{12}. Hence, we obtain for every $\ell \in \left(k-\frac{1}{2p}, k\right]$ the following bound for $F_{\ell}$:
\begin{equation}\label{bound_p_const}
    \begin{split}
         F_{\ell} & \leq \binom{n}{k} \binom{k}{\ell} \binom{n-k}{k-\ell} p^{2(k-1)} (1-p)^{2 \binom{k}{2} - \binom{\ell}{2} -2(k-1)} k^{k-2} \ell \max_{r \in\{0,\ldots,\ell-1\}}f(k,\ell,r)  \frac{(1-p)^{r}}{p^r}\\
         & = \binom{n}{k} \binom{k}{\ell} \binom{n-k}{k-\ell} p^{2(k-1)} (1-p)^{2 \binom{k}{2} - \binom{\ell}{2} -2(k-1)} k^{k-2} \ell \max_{r \in\{0,\ldots,\ell-1\}}g(k,\ell,r),
    \end{split}
\end{equation}
where $g(k,\ell,r) = f(k, \ell, r) \left( \frac{1-p}{p}\right)^r$. Let us find a suitable upper bound for $g(k,\ell,r)$ following the same method that in \cite{12}. We will compute the derivative of \( g(k, \ell, r) \) with respect to \( r \) and determine a point \( \hat{r} \) at which the function attains its maximum for any given \( \ell \) within the considered interval.

Let $\hat{r}:=\ell \left(1-\frac{p(k-\ell)}{(1-p)e} \right)$. We get that for any $\ell \in \left(k-\frac{1}{2p}, k-1\right]$, it holds that 
\begin{equation*}
       \frac{\partial}{\partial r} \ln g(k,\ell,r)=\ln \underbrace{\left( \frac{2(1-p)}{(k-\ell)p}\right)}_{>1}I[r< \ell/2]+\ln \underbrace{\left( \frac{9(1-p)}{8(k-\ell)p}\right)}_{>1} I[\ell/2 \leq r < \ell (1-1/e)]
\end{equation*}

\begin{equation*}
    +\ln\underbrace{\left(\frac{e(\ell-r)(1-p)}{\ell(k-\ell)p} \right)}_{>1 \text{ for } r<\hat{r} \text{ ,and } <1 \text{ for } r>\hat{r}} I[r\geq \ell (1-1/e)],
\end{equation*}
and due to (\ref{equalities}), (\ref{left_1}) and (\ref{right_1}), we thus get that \( g(k,\ell,r) \) increases in \( r \in [0,\hat{r}] \) and decreases in \( r \in [\hat{r},\ell-1] \). Therefore, \( g(k,\ell,r) \) is maximal at \( g(k,\ell,\hat{r}) \).\\

Thus, given that from \cite{12},
$$
f(k,\ell,r)=\left(\frac{\ell}{\ell-r} \right)^{\ell-r}(\ell+1)^{k-\ell-1}(k-\ell)^{k-r-2}\text{ for }r \in \left[\ell(1-1/e), \ell-1\right),
$$
we get that for any \( \ell \in \left(k-\frac{1}{2p}, k\right] \), it holds that

\begin{equation*}
\begin{split}
    g(k,\ell,r)\leq g(k,\ell,\hat{r}) & =\left(\frac{1-p}{p}\right)^{\hat{r}}\left(\frac{\ell}{\ell-\hat{r}} \right)^{\ell-\hat{r}}(\ell+1)^{k-\ell-1}(k-\ell)^{k-2-\hat{r}}\\
    & =\left(\frac{1-p}{p}\right)^{\ell \left(1-\frac{p(k-\ell)}{(1-p)e} \right)}\left(\frac{e(1-p)}{p(k-\ell)} \right)^{\ell\frac{p(k-\ell)}{(1-p)e}}(\ell+1)^{k-\ell-1}(k-\ell)^{k-2-\ell+\ell\frac{p(k-\ell)}{(1-p)e}}\\
    & = \left(\frac{1-p}{p}\right)^{\ell}(\ell+1)^{k-\ell-1}(k-\ell)^{k-\ell-2}e^{\frac{p}{(1-p)}(k-\ell)\frac{\ell}{e}}.
\end{split}
\end{equation*}
In this way, we obtain the following upper bound for (\ref{ecuacion_4_1}) within the considered interval: 
\begin{equation}
\begin{split}
    \displaystyle \sum_{\ell = \lceil k-\frac{1}{2p}\rceil }^{k} \frac{F_{\ell}}{({\sf E} X_k)^2}  & \leq  \displaystyle \sum_{\ell = \lceil k-\frac{1}{2p}\rceil}^{k} \frac{\binom{k}{\ell} \binom{n-k}{k-\ell}}{\binom{n}{k} k^{k-2}} (1-p)^{-\binom{\ell}{2}} \ell \left(\frac{1-p}{p}\right)^{\ell}(\ell+1)^{k-\ell-1}(k-\ell)^{k-\ell-2} e ^{\frac{\ell(k-\ell)p}{e(1-p)}} \\
    & = \displaystyle \sum_{\ell = \lceil k-\frac{1}{2p}\rceil}^{k}  \hat{I}_{\ell},
\end{split} \label{eq_final}
\end{equation}
where 
$$
\hat{I}_{\ell}=\frac{\binom{k}{\ell} \binom{n-k}{k-\ell}}{\binom{n}{k} k^{k-2}} (1-p)^{-\binom{\ell}{2}} \ell \left(\frac{1-p}{p}\right)^{\ell}(\ell+1)^{k-\ell-1}(k-\ell)^{k-\ell-2} e ^{\frac{\ell(k-\ell)p}{e(1-p)}}.
$$
Let us show that $\hat{I}_{\ell}$ increases in $\ell \in \left(k-\frac{1}{2p}, k-1\right]$. For this, we will show that $\frac{\hat{I}_{\ell+1}}{\hat{I}_{\ell}}>1$ within this interval and for $p\geq n^{-\theta}$, $\theta< \frac{e-2}{3e-2} $. For simplicity of notation, let us write $k-\beta$ instead of $\ell$, where $\beta \in [1,\frac{1}{2p}]$. We thus get
\begin{equation*}
    \begin{split}
        \frac{\hat{I}_{\ell+1}}{\hat{I}_{\ell}}& =  \frac{(k-\ell)^2(1-p)^{1-\ell}}{\ell(\ell+2)(n-2k+\ell+1)(k-\ell-1)p}\left(\frac{\ell+2}{\ell+1} \right)^{k-\ell-1} \left(\frac{k-\ell-1}{k-\ell} \right)^{k-\ell-2} e^{\frac{p(k-2\ell-1)}{(1-p)e}} \\
        & > \frac{\beta}{(\ell+2)^2n}\frac{(1-p)^{1-\ell}}{p}\left( \frac{\beta-1}{\beta}\right)^{\beta-2}e^{\frac{p}{e(1-p)}(-k+2\beta-1)}\\
        & = e^{\ln \beta-2\ln(\ell+2)-\ln n +(1-\ell)\ln(1-p)-\ln p +(\beta-2)\ln(\frac{\beta-1}{\beta})+\frac{p}{(1-p)e}(-k+2\beta-1)}\\
        & = e^{\ln \beta-2\ln(k-\beta+2)-\ln n +(1-k+\beta)\ln(1-p)-\ln p +(\beta-2)\ln(\frac{\beta-1}{\beta})+\frac{p}{(1-p)e}(-k+2\beta-1)}\\
        & > e^{\ln \beta-2\ln k-\ln n -k\ln(1-p)+(1+\beta)\ln(1-p)-\ln p +(\beta-2)\ln(\frac{\beta-1}{\beta})-\frac{pk}{(1-p)e}+\frac{p}{(1-p)e}(2\beta-1)}\\
        & > e^{\ln \left(\frac{1}{k^2 np} \right)+kp-\frac{pk}{(1-p)e}+\ln 2 +(1+\beta)\ln(1-p) +(\beta-2)\ln(\frac{\beta-1}{\beta})+\frac{p}{(1-p)e}(2\beta-1)}\\
        & = e^{\ln \left(\frac{1}{k^2 np} \right)+kp\left[1-\frac{1}{(1-p)e}\right]+\ln 2 +(1+\beta)\ln(1-p) +(\beta-2)\ln(\frac{\beta-1}{\beta})+\frac{p}{(1-p)e}(2\beta-1)}
\end{split}
\end{equation*}
Substituting the value of $k$ into the expression, we obtain
\begin{equation*}
    \begin{split}
        \frac{\hat{I}_{\ell+1}}{\hat{I}_{\ell}}& >  e^{\ln \left(\frac{p}{16 n (\ln (np))^2} \right)+(2\ln(npe)+o(1))\left[1-\frac{1}{(1-p)e}\right]+\ln 2 +(1+\beta)\ln(1-p) +(\beta-2)\ln(\frac{\beta-1}{\beta})+\frac{p}{(1-p)e}(2\beta-1)}\\
        & = e^{\ln \left(\frac{p}{n} \right)+2\ln(npe)\left[1-\frac{1}{(1-p)e}\right]-\ln(16(\ln(np))^2)+\ln 2 +(1+\beta)\ln(1-p) +(\beta-2)\ln(\frac{\beta-1}{\beta})+\frac{p}{(1-p)e}(2\beta-1)}\\
        & > e^{\ln n \left[-(\theta+1)+2(1-\theta)\left[1-\frac{1}{e}\right] \right]-\ln(16(\ln(np))^2)+\ln 2 +(1+\beta)\ln(1-p) +(\beta-2)\ln(\frac{\beta-1}{\beta})+\frac{p}{(1-p)e}(2\beta-1)}\\
        & \to \infty, \ \text{ for } \theta<\frac{e-2}{3e-2}.
    \end{split}
    \label{ec_il}
\end{equation*}
As a result, if $p\geq n^{-\theta}$, $\theta< \frac{e-2}{3e-2} $, then $\frac{\hat{I}_{\ell+1}}{\hat{I}_{\ell}}>1$ for all $\ell \in \left(k-\frac{1}{2p}, k-1\right]$.

Therefore, we can bound the right-hand side of (\ref{eq_final}) by $k\cdot \hat{I}_{\ell}$ at $\ell=k-1$. This results in the following bound, obtained using Stirling’s formula:
\begin{equation*}
\begin{split}
    \displaystyle \sum_{\ell = \lceil k-\frac{1}{2p}\rceil}^{k-1} \frac{F_{\ell}}{({\sf E} X_k)^2}   \leq k \cdot \hat{I}_{\ell}{\bigg |}_{\ell=k-1}
    & \leq k \cdot \frac{k(n-k)}{\binom{n}{k}k^{k-2}}(1-p)^{-\binom{k-1}{2}}(k-1)\left(\frac{1-p}{p} \right)^{k-1}e^{\frac{(k-1)p}{e(1-p)}} \\
    & \leq  \frac{(n-k)}{\binom{n}{k}k^{k-5}}(1-p)^{-\binom{k-1}{2}}\left(\frac{1-p}{p} \right)^{k-1}e^{\frac{(k-1)p}{e(1-p)}}\\
    & \sim \sqrt{2\pi} \frac{(n-k)k^{\frac{11}{2}}}{n^k e^k }(1-p)^{-\binom{k-1}{2}}\left(\frac{1-p}{p} \right)^{k-1}e^{\frac{(k-1)p}{e(1-p)}},
\end{split} 
\end{equation*}
and applying the exponential function to the entire expression gives
\begin{equation}
\begin{split}
    \displaystyle \sum_{\ell = \lceil k-\frac{1}{2p}\rceil}^{k-1} \frac{F_{\ell}}{({\sf E} X_k)^2}  & \leq e^{\frac{1}{2}\ln 2\pi +\ln (n-k)+\frac{11}{2}\ln k -k\ln n -k +(k-1-\binom{k-1}{2})\ln (1-p) +(1-k)\ln p +\frac{(k-1)p}{e(1-p)}}\\
    &= e^{ -k \ln n -\frac{k^2}{2}\ln (1-p)-k\ln p-k+\ln (n-k)+\frac{11}{2}\ln k -kp -\frac{3}{2}kp +\ln p +\frac{kp}{e(1-p)}+O(1)+O(kp^2)}\\
    &= e^{ -k \ln (np) +\frac{k^2}{2}p+\frac{k^2p^2}{4} -k+\ln n +\frac{11}{2}\ln k +kp\left[-\frac{5}{2}+\frac{1}{e(1-p)}\right] +\ln p+O(1)+O(kp^2)+O(k^2p^3)}\\
    &= e^{ k\left[- \ln (np) +\frac{k}{2}p+\frac{kp^2}{4} -1+ \frac{1}{p}\left[-\frac{5}{2}+\frac{1}{e(1-p)}\right]\right]+\ln (np) +\frac{11}{2}\ln k +O(1)+O(kp^2)+O(k^2p^3)}\\
    &= e^{ k\left[- \ln (np) + \ln(npe)+p\frac{\ln(npe)}{2} -1+ \frac{1}{p}\left[-\frac{5}{2}+\frac{1}{e(1-p)}\right]\right]+\ln (np^{-9/2}) +O(\ln \ln n)+O(kp^2)+O(k^2p^3)}\\
    &= e^{k \left[\frac{p}{2}\ln(np)  +\frac{1}{p}\left(-\frac{5}{2}+\frac{1}{e}\right)+o(1)\right]+\ln (np^{-9/2}) +O(\ln \ln n)+O(kp^2)+O(k^2p^3)}.
\end{split} 
\label{lastone}
\end{equation}
Since for any $p$ such that $n^{-\theta} \leq p<\frac{1}{2\ln n}$, $\theta< \frac{e-2}{3e-2}$, the dominant asymptotics in the right-hand side of (\ref{lastone}) will be determined by $e^{\frac{k}{p}\left(-\frac{5}{2}+\frac{1}{e}\right)}$, then the fact that $\frac{k}{p}\left(-\frac{5}{2}+\frac{1}{e}\right) \to -\infty$, implies that the right-hand side of (\ref{eq_final}) is o(1). \\

Finally, using all the previous results obtained in (\ref{eq_final_1}), (\ref{eq_final_2}), and  (\ref{eq_final_3}), we obtain that for any \( p \) such that $n^{-\frac{e-2}{3e-2}+\varepsilon}\leq p<\frac{1}{2\ln n}$, $\varepsilon>0$, it holds that
\begin{equation*}
    {\sf P}(X_k=0)\leq \frac{\mathrm{Var}X_k}{({\sf E}X_k)^2}\leq \displaystyle \sum_{\ell =2}^{k} \frac{F_{\ell}}{({\sf E}X_k)^2}=o(1).
\end{equation*}
This concludes the proof of Theorem~\ref{th:new} for  $n^{-\frac{e-2}{3e-2}+\varepsilon} \leq p < \frac{1}{2\ln n}$.


\section{Proof of Theorem~\ref{th:new} for $\frac{1}{2\ln n} \leq p = o(1)$}\label{p_large}

In this case, we adopt the techniques from \cite{12}, originally applied to constant $p$, to show that the right-hand side of (\ref{ecuacion_4_1}) is $o(1)$. For the sake of completeness, we will present the full adapted version of the proof here.

Similarly to what we did in Sections~\ref{sc:part2} and~\ref{sc:part4}, we bound $N(k, \ell, r) \leq k^{k-2} f(k, \ell, r)$. Consequently, we obtain the bound for $F_{\ell}$ given in (\ref{bound_p_const}).  

Set $\ell_1=2\frac{\ln n}{\ln[1/(1-p)]}-16\frac{\ln \ln n}{\ln[1/(1-p)]}$. For $\ell\leq \ell_1$, we will use the trivial bound $f(k,\ell,r)=k^{k-2}$. For such $\ell$,
\begin{equation*}
    \frac{F_{\ell}}{(EX_k)^2} \leq \frac{\binom{k}{\ell}\binom{n-k}{k-\ell}(1-p)^{-\binom{\ell}{2}}\max_{r\in\{0,\ldots,\ell-1\}} ((1-p)/p)^{r}}{\binom{n}{k}}.
\end{equation*}
Since $1-p>p$, and due to (\ref{bound_11}) we thus get
\begin{equation*}
\begin{split}
       \frac{F_{\ell}}{(EX_k)^2} \leq & \frac{\binom{k}{\ell}\binom{n-k}{k-\ell}(1-p)^{-\binom{\ell}{2}} ((1-p)/p)^{\ell}}{\binom{n}{k}}\leq \sqrt{2\pi k}\left(\frac{k^2 e (1-p)}{\ell np}(1-p)^{-\ell/2}\right)^{\ell}\\
       & =e^{\frac{1}{2}\ln(2\pi k)+\ell \left[2\ln k+1+\ln(1-p)-\frac{\ell}{2}\ln(1-p)-\ln \ell-\ln (pn) \right]}\\
       & = e^{\frac{1}{2}\ln(2\pi k)+\ell \left[2\ln k+1-p+O(p^2)+\frac{\ell}{2}\ln(1/(1-p))-\ln \ell-\ln (pn) \right]}\\
       & \leq e^{\frac{1}{2}\ln(2\pi k)+\ell \left[2\ln k+1-8\ln \ln n-\ln \ell-\ln p +o(1)\right]}\\
       & \leq e^{\frac{1}{2}\ln(2\pi k)+\ell \left[2\ln 2+2\ln \ln n-2\ln p+1-8\ln \ln n-\ln 2-\ln p +o(1)\right]}\\
       & = e^{\frac{1}{2}\ln(2\pi k)+\ell \left[\ln \left( \frac{2e}{[(\ln n)^2 p]^3}\right)+o(1)\right]}.
\end{split} 
\end{equation*}
Therefore, having $p\geq \frac{1}{2\ln n}$, we obtain
\begin{equation*}
    \displaystyle \sum_{\ell =2}^{\ell_1} \frac{F_{\ell}}{({\sf E}X_k)^2}\leq k e^{\frac{1}{2}\ln(2\pi k)+2 \left[\ln \left( \frac{2e}{[(\ln n)^2 p]^3}\right)+o(1)\right]}=o(1). 
\end{equation*}
Let us now consider the case when $\ell \in \left(\ell_1, k-\frac{2(1-p)}{p} \right]$. Here we bound $N(k, \ell, r)\leq k^{k-2}f(k,\ell,r)$ and by following the same arguments as in \cite{12}, we obtain the bound for $F_{\ell}$ given in (\ref{bound_2}), that we denoted by $\hat{F}_{\ell}$. Moreover, since $\frac{\partial}{\partial \ell} \ln \left[\frac{\hat{F}_{\ell+1}}{\hat{F}_{\ell}}\right]<0$, then $\frac{\hat{F}_{\ell+1}}{\hat{F}_{\ell}}$ dicreases with $\ell$ in the range.
Let us show that 
\begin{equation*}
    \max_{\ell \in\{\ell_1+1,\ldots,k-\frac{3(1-p)}{p}\}}\hat{F}_{\ell}=\hat{F}_{k-\frac{3(1-p)}{p}}.
\end{equation*}
For this, it suffices to show $\frac{\hat{F}_{\ell+1}}{\hat{F}_{\ell}}\geq 1$ for $\ell=k-\frac{3(1-p)}{p}$. For the sake of notation, let $\ell_2=k-\frac{3(1-p)}{p}$.

In this way, since we have 
\begin{itemize}
\item $(1-p)^{-\ell_2}= \left( \frac{1}{1-p}\right)^{\frac{\ln \left(\frac{n^2p^2}{e} \right)}{\ln(1/[1-p])}+6+o(1)}\geq \frac{n^2p^2}{e^2}$;
\item $\left(\frac{k-\ell_2-1}{k-\ell_2} \right)^{k-2} =e^{(k-2)\ln(1-\frac{1}{k-\ell_2})}\geq e^{-\frac{kp}{3(1-p)}-O\left(kp^2 \right)}=e^{-\frac{2}{3}\ln (npe)-O\left(kp^2 \right)}=\frac{e^{-O\left(kp^2 \right)}}{(npe)^{2/3}}$;
\item $\frac{\hat{F}_{\ell+1}}{\hat{F}_{\ell}} =\frac{(k-\ell)^2(1-p)^{-\ell}}{(\ell+1)^2(n-2k+\ell+1)}\left(\frac{k-\ell-1}{k-\ell}\right)^{k-2}\left(\frac{\ell+2}{\ell+1} \right)^{k-\ell-2}$,
\end{itemize} 
we get
\begin{equation*}
    \begin{split}
         \frac{\hat{F}_{\ell_2+1}}{\hat{F}_{\ell_2}} & \geq \frac{9(1-p)^2}{p^2}\frac{1}{\ell_2^2}\frac{1}{n}\frac{n^2p^2}{e^2}(npe)^{-2/3}e^{-O\left(kp^2 \right)}=\frac{(1-p)^2p^{-2/3}e^{-O\left(kp^2 \right)}}{\ell_2^2}{ n^{1/3}} \to \infty.
    \end{split} 
\end{equation*}

Let us now find the asymptotics of $\hat{F}_{\ell}$ for any $\ell = k - \frac{c (1 - p)}{p}$, where $c \geq 2$. Due to (\ref{bound_12})~we~get
\begin{equation*}
    \begin{split}
        \hat{F}_{\ell}& =\frac{\binom{k}{\ell} \binom{n-k}{k-\ell} (1-p)^{-\binom{\ell}{2}} (k-\ell)^{k-2} (\ell+1)^{k-\ell-1}}{\binom{n}{k} k^{k-3}}  \\
        & \leq \left( \frac{k^2}{\ell ne} \right)^{\ell}\left(\frac{k}{k-\ell}\right)^{2(k-\ell)} \frac{(k-\ell)^{k-2}(\ell+1)^{k-\ell-1}(1-p)^{-\binom{\ell}{2}}}{k^{k-3}}\\
        & =\frac{k^{k+3}(k-\ell)^{2\ell-k-2}(\ell+1)^{k-\ell-1}(1-p)^{-\binom{\ell}{2}}}{(\ell e n)^{\ell}}\\
        & \leq \frac{k^{k+3}}{(\ell e n)^{\ell}}\left(\frac{c}{p} \right)^{k-\frac{2c}{p}+2c-2}(\ell+1)^{\frac{c}{p}-1}(1-p)^{-\binom{\ell}{2}} =e^{g(\ell)},
    \end{split}
\end{equation*}
where
\begin{equation*}
    \begin{split}
         g(\ell)& =-\ell \ln n-\frac{\ell^2}{2}\ln(1-p)-\ell \ln \ell+(k+3)\ln k+k\ln \left(\frac{c}{p} \right)-\ell+O\left(\frac{\ln \ln n}{p}\right)\\
        & = -\ell \left[\ln n+\frac{\ell}{2}\ln (1-p) \right]-\left(k-\frac{c(1-p)}{p}\right)\ln \ell +k\ln k-k\ln p+k\ln c-\ell+O\left(\frac{\ln \ln n}{p}\right)\\
         & = -\ell \left[\ln n-\frac{k-[c(1-p)/p]}{2}\ln \left( \frac{1}{1-p}\right) \right]+k \ln (k/\ell)-k\ln p+k\ln c-\ell+O\left(\frac{\ln \ln n}{p}\right)\\
          & = -\ell \left[-\ln(pe)+\frac{c}{2p}\ln\left( \frac{1}{1-p}\right)+O(p) \right]+o(k)-k\ln p+k\ln c-\ell+O\left(\frac{\ln \ln n}{p}\right).
\end{split}
\end{equation*}
Furthermore, by simplifying the expression, we obtain 
\begin{equation*}
    \begin{split}
           g(\ell) & = \ell\ln(pe)-\ell\frac{c}{2p}\ln\left( \frac{1}{1-p}\right)+O(\ell p) +o(k)-k\ln p+k\ln c-\ell+O\left(\frac{\ln \ln n}{p}\right)\\
           & = k\ln p+k -\ell\frac{c}{2p}\ln\left( \frac{1}{1-p}\right)+O(\ell p) +o(k)-k\ln p+k\ln c-\ell+O\left(\frac{\ln \ln n}{p}\right)\\
            & = k -\ell\frac{c}{2p}\ln\left( \frac{1}{1-p}\right)+O(\ell p) +o(k)+k\ln c-k+\frac{c(1-p)}{p}+O\left(\frac{\ln \ln n}{p}\right)\\
            & = -k\frac{c}{2p}\ln\left( \frac{1}{1-p}\right)+k\ln c+O(\ell p) +o(k)+O\left(\frac{\ln \ln n}{p}\right)\\
            & = k\underbrace{\left[-\frac{c}{2p}\ln\left( \frac{1}{1-p}\right)+\ln c\right]}_{\boldsymbol{<0}\text{ for any } c>0}+O(\ell p) +o(k)+O\left(\frac{\ln \ln n}{p}\right).
    \end{split}
\end{equation*}
Therefore, we obtain
\begin{equation*}
    \displaystyle \sum_{\ell =\ell_1+1}^{k-\frac{2(1-p)}{p}} \frac{F_{\ell}}{({\sf E}X_k)^2}\leq \displaystyle \sum_{\ell =\ell_1+1}^{\ell_2}\hat{F_{\ell}}  + \displaystyle \sum_{\ell =\ell_2+1}^{k-\frac{2(1-p)}{p}}\hat{F_{\ell}}   \leq 2ke^{g(\ell)}=o(1).
\end{equation*}

Finally, we can move to the case when $\ell \in \left( k-\frac{2(1-p)}{p},k-1 \right]$. Set $\ell=k-s$, where $s \in \left[1, \frac{2(1-p)}{p}\right)$. Similarly to~\cite{12}, we get that for some constant $A>0$,
\begin{equation}\label{thm_1_last}
    \frac{F_{\ell}}{({\sf E}X_k)^2}\leq \frac{A \binom{k}{s}\binom{n-k}{s}(1-p)^{sk}k^s}{{\sf E}X_k} \max_{r\in \{ 0,1,\ldots,k-s-1\}}f_0(k,r)s^{k-r}\left(\frac{p}{1-p} \right)^{k-r},
\end{equation}
where 
\begin{equation*}
    f_0(k,r):=\left(\frac{\ell}{\ell-r} \right)^{k-r} I(r\geq \ell(1-1/e))+(4/3)^k (9/8)^r I(\ell/2 \leq r <\ell(1-1/e))+2^rI(r<\ell/2).
\end{equation*}
Furthermore, applying Stirling’s approximation, we observe that
\begin{equation*}
    \begin{split}
         \binom{k}{s}\binom{n-k}{s}(1-p)^{sk}k^s & \leq \frac{1}{(s!)^2}k^{2s}n^{-s}(pe)^{-2s}\left( \frac{1}{1-p}\right)^{-3s+o(1)}\\
         & \leq n^{-s}e^{-(2s+1)\ln s+2s\ln k -2s \ln(pe)+O(1/p)}\\
          & \leq e^{-s \ln n+O(s \ln \ln n)},
    \end{split}
\end{equation*}
As a consequence, (\ref{thm_1_last}) leads to
\begin{align*}
    \frac{F_{\ell}}{({\sf E}X_k)^2}\leq \frac{\exp \left[-s \ln n+\displaystyle\max_{r\in \{ 0,1,\ldots,k-s-1\}}\ln f_1(k,r)+O(s \ln \ln n) \right]}{{\sf E}X_k},\quad\text{ where} \\
    f_1(k,r)=f_0(k,r)s^{k-r}\left(\frac{p}{1-p} \right)^{k-r}.
\end{align*}

The next step is to investigate how the value of $r$ that maximizes $f_1(k, r)$ depends on $\ell$, noting that the maximum may be attained at different values of $r$ for different intervals of $\ell$.

To facilitate the analysis, we partition the values of $\ell \in \left(k-\frac{2(1-p)}{p},k-1 \right]$ into two distinct zones as follows:
\begin{equation*}
   \underbrace{\left(k-\frac{2(1-p)}{p}, k-\frac{9}{8}\frac{(1-p)}{p}\right)}_{\text{interval 1}}\sqcup \underbrace{\left[k-\frac{9}{8}\frac{(1-p)}{p}, k-1\right)}_{\text{interval 2}}.
\end{equation*}
In this way, for $r<\ell/2$, since $\frac{\partial \ln f_1(k,r)}{\partial r}=\ln \frac{2(1-p)}{ps}$, we get that $f_1(k,r)$ increases in both intervals. On the other hand, for $r\in [\ell/2,\ell(1-1/e))$, since $\frac{\partial \ln f_1(k,r)}{\partial r}=\ln \frac{9}{8}\frac{(1-p)}{ps}$, we get that $f_1(k,r)$ decreases in interval 1 and increases in interval 2. Finally, for $r\in [\ell(1-1/e),\ell-1]$, set $r=\ell-x$, where $1\leq x\leq \ell/e$, and we thus get that $\ln f_1(k,r)$ is given by
\begin{equation*}
   (k-r)\ln \left( \frac{\ell}{(\ell-r)}\frac{ps}{1-p}\right)=(s+x)\ln \left( \frac{\ell }{x}\frac{ps}{1-p}\right)=\underbrace{s\ln \left( \frac{\ell}{x}\frac{ps}{1-p}\right)}_{o(s \ln n)}+x\ln \left( \frac{\ell}{x}\frac{ps}{1-p}\right),
\end{equation*}
where only the asymptotics of $x\ln \left( \frac{\ell}{x}\frac{ps}{1-p}\right)$ will be relevant. Then, since
\begin{equation*}
    \frac{\partial}{\partial x}\left[x\ln \left( \frac{\ell}{x}\frac{ps}{1-p}\right)\right]=\ln\left( \frac{\ell}{x} \frac{ps}{(1-p)}\right)-1,
\end{equation*}
in both intervals, the function $f_1(k, r)$ first increases and then decreases, attaining its maximum at $r\sim \ell\left(1-\frac{ps}{(1-p)e}\right)$.

Notice that $f_1(k,r)$ has two points of discontinuity: $r_1=\ell/2$, $r_2=\ell(1-1/e)$. It is straightforward to verify that $f_1(k,r_1)>f_1(k,r_1-0)$ and $f_1(k,r_2)>f_1(k,r_2-0)$.

In summary, we have shown that
\begin{itemize}
    \item In interval 2, $f_1(k,r)$ achieves its maximum at $r\sim \ell\left(1-\frac{ps}{(1-p)e}\right)$.
    \item In interval 1, $f_1(k,r)$ achieves its maximum either at $r=\ell/2$, or at $r\sim \ell\left(1-\frac{ps}{(1-p)e}\right)$.
\end{itemize}
We now proceed to determine the asymptotics of the sum of $\frac{F_{\ell}}{({\sf E}X_k)^2}$ for all $\ell$ in interval 2.
\begin{equation}\label{thm_last_2}
    \begin{split}
       &\displaystyle \sum_{\ell=k-\frac{9}{8}\frac{(1-p)}{p}}^{k-1} \frac{F_{\ell}}{({\sf E}X_k)^2}\leq \displaystyle \sum_{\ell=k-\frac{9}{8}\frac{(1-p)}{p}}^{k-1}  \frac{\exp \left[-s \ln n+\displaystyle\max_{r\in \{ 0,1,\ldots,k-s-1\}}\ln f_1(k,r)+O(s \ln \ln n) \right]}{{\sf E}X_k}\\
        &=\displaystyle \sum_{\ell=k-\frac{9}{8}\frac{(1-p)}{p}}^{k-1} \frac{\exp \left[-s \ln n+\displaystyle\max_{r\in \{ 0,1,\ldots,k-s-1\}}(k-r)\ln\left(\frac{\ell}{\ell-r}\frac{ps}{(1-p)} \right)+O(s \ln \ln n) \right]}{{\sf E}X_k}\\
        &=\displaystyle \sum_{\ell=k-\frac{9}{8}\frac{(1-p)}{p}}^{k-1} \frac{\exp \left[-s \ln n+\left[\frac{\ell ps}{(1-p)e} -s\right]+O(s \ln \ln n) \right]}{{\sf E}X_k}\\
        &= \displaystyle \sum_{\ell=k-\frac{9}{8}\frac{(1-p)}{p}}^{k-1}\frac{\exp \left[-s \ln n\left(1+\frac{1}{\ln n}-\frac{\ell}{\ln n}\frac{p}{(1-p)e} \right)+O(s \ln \ln n) \right]}{{\sf E}X_k}\\
         &\leq k\frac{\exp \left[-s \ln n\left(1-\frac{2}{(1-p)e}+o(1) \right)+O(s \ln \ln n) \right]}{{\sf E}X_k}=o(1).
    \end{split}
\end{equation}
On the other hand, in interval 2, if $f_1(k, r)$ attains its maximum at $r \sim \ell\left(1 - \frac{ps}{(1 - p)e} \right)$, then, by an argument similar to that used in (\ref{thm_last_2}), we obtain that
\begin{equation*}
    \displaystyle \sum_{\ell>k-\frac{2(1-p)}{p}}^{k-\frac{9}{8}\frac{(1-p)}{p}} \frac{F_{\ell}}{({\sf E}X_k)^2}=o(1).
\end{equation*}
Finally, if in interval 2, $f_1(k, r)$ attains its maximum at $r=\ell/2$, then we get the following: 
\begin{multline}\label{thm_last_3_2}
\sum_{\ell>k-\frac{2(1-p)}{p}}^{k-\frac{9}{8}\frac{(1-p)}{p}} \frac{F_{\ell}}{({\sf E}X_k)^2} \leq \displaystyle \sum_{\ell=k-\frac{9}{8}\frac{(1-p)}{p}}^{k-1}  \frac{\exp \left[-s \ln n+\displaystyle\max_{r\in \{ 0,1,\ldots,k-s-1\}}\ln f_1(k,r)+O(s \ln \ln n) \right]}{{\sf E}X_k}\\
=\sum_{\ell=k-\frac{9}{8}\frac{(1-p)}{p}}^{k-1}  \frac{\exp \left[-s \ln n+\displaystyle\max_{r\in \{ 0,1,\ldots,k-s-1\}}\ln \left((4/3)^k(9/8)^r \left(\frac{ps}{1-p}\right)^{k-r}\right)+O(s \ln \ln n) \right]}{{\sf E}X_k}\\
=\sum_{\ell=k-\frac{9}{8}\frac{(1-p)}{p}}^{k-1}  \frac{\exp \left[-s \ln n+(k-r)\ln \left(\frac{4}{3}\frac{ps}{(1-p)}\right)-r\ln (2/3)+O(s \ln \ln n) \right]}{{\sf E}X_k},
\end{multline}
and replacing $r$ by $\ell/2$ we get that the right-hand side of \eqref{thm_last_3_2} is bounded by
\begin{multline}\label{thm_last_3}
\sum_{\ell=k-\frac{9}{8}\frac{(1-p)}{p}}^{k-1}  \frac{\exp \left[-s \ln n+(k-\ell/2)\ln \left(\frac{4}{3}\frac{ps}{(1-p)}\right)-\frac{\ell}{2}\ln (2/3)+O(s \ln \ln n) \right]}{{\sf E}X_k}\\
         \leq k  \frac{\exp \left[-s \ln n+\frac{k}{2}\ln \left(\frac{4}{3}\frac{ps}{(1-p)}\frac{3}{2}\right)+O(s \ln \ln n) \right]}{{\sf E}X_k}\\
         =  \frac{\exp \left[-\ln n\left[s-\frac{1}{p}\ln\left(\frac{2ps}{1-p} \right)+o(1) \right]
          \right]}{{\sf E}X_k},
\end{multline}
where $s-\frac{1}{p}\ln\left(\frac{2ps}{1-p} \right)$ is minimal at $s=\frac{9}{8}\frac{(1-p)}{p}$, since $\frac{\partial}{\partial s}\left[s-\frac{1}{p}\ln\left(\frac{2ps}{1-p} \right) \right]>0$. Therefore, (\ref{thm_last_3}) gives

\begin{equation*}
    \begin{split}
       \displaystyle \sum_{\ell>k-\frac{2(1-p)}{p}}^{k-\frac{9}{8}\frac{(1-p)}{p}} \frac{F_{\ell}}{({\sf E}X_k)^2}\leq  \frac{\exp \left[-\ln n\left[s-\frac{1}{p}\ln\left(\frac{2ps}{1-p} \right)+o(1) \right]
          \right]}{{\sf E}X_k}\leq \frac{\exp \left[-\ln n\left[\frac{9}{8}-\ln (9/4)\right]
          \right]}{{\sf E}X_k}=o(1).
    \end{split}
\end{equation*}
This concludes the proof of Theorem \ref{th:new} for  $\frac{1}{2\ln n} \leq p = o(1)$.

\paragraph{Acknowledgments} The author is deeply grateful to Maksim Zhukovskii for his constant support and valuable feedback on the presentation of the paper.

\end{document}